\documentclass{article}
\usepackage{iclr2027_conference,times}

\usepackage{amsmath,amsfonts,bm}

\def\eqref#1{equation~\ref{#1}}

\def\1{\bm{1}}

\DeclareMathAlphabet{\mathsfit}{\encodingdefault}{\sfdefault}{m}{sl}
\SetMathAlphabet{\mathsfit}{bold}{\encodingdefault}{\sfdefault}{bx}{n}

\newcommand{\E}{\mathbb{E}}

\usepackage{hyperref}
\usepackage{url}
\hypersetup{hidelinks,hypertexnames=false}

\title{Optimal Momentum Methods for Stochastic\\ Multilevel Compositional Optimization}

\author{
{Wei Jiang}$^{1}$ \quad
{Rui Yan}$^{1}$ \quad
{Sifan Yang}$^{2}$ \quad
{Yuanyu Wan}$^{3}$ \quad
{Lijun Zhang}$^{2}$ \quad
{Zechao Li}$^{1}$ \\[0.5em]
$^{1}$School of Computer Science and Engineering, \\
Nanjing University of Science and Technology, Nanjing, China \\
$^{2}$School of Artificial Intelligence, Nanjing University, Nanjing, China \\
$^{3}$School of Software Technology, Zhejiang University, Ningbo, China
}
\usepackage{amsmath}
\usepackage{amssymb}
\usepackage{mathtools}
\usepackage{amsthm}
\usepackage{subfigure}
\theoremstyle{plain}
\newtheorem{theorem}{Theorem}[section]

\newtheorem{lemma}[theorem]{Lemma}

\theoremstyle{definition}

\newtheorem{ass}[theorem]{Assumption}
\theoremstyle{remark}

\usepackage{bbding}
\usepackage{microtype}
\usepackage{graphicx}
\usepackage{booktabs}
\usepackage{paralist}
\usepackage{multirow}
\RequirePackage{algorithm}
\RequirePackage{algorithmic}

\newcommand{\Norm}[1]{\left\|#1\right\|}

\def \cite {\citep}
\def \rat {0.27}
\def \rata {0.325}
\def \E {\mathbb{E}}

\def \r {\mathbf{r}}

\def \u {\mathbf{u}}
\def \z {\mathbf{v}}

\def \x {\mathbf{w}}

\def \z {\mathbf{z}}
\def \x {\mathbf{x}}
\def \y {\mathbf{y}}

\def \LL {L} 

\usepackage{textcomp}
\iclrfinalcopy
\begin{document}

\maketitle

\begin{abstract}
This paper investigates stochastic multi-level optimization where the objective is a nested composition of several smooth non-convex functions.
We assume that only stochastic estimates of the gradient and function values for each level are accessible.
Consequently, obtaining an accurate estimate of the overall gradient is challenging due to the nested structure. 
To address this, we employ a momentum-based estimator with mini-batches to track the function values of each level, which are subsequently used to construct momentum gradient estimators.
We establish an optimal sample complexity of $\mathcal{O}(\epsilon^{-4})$ for finding an $\epsilon$-stationary point, avoiding the stronger average smoothness assumption commonly relied upon in prior literature.
Furthermore, by employing a normalization technique, we attain the same rate without requiring problem-dependent constants to set hyperparameters.
To achieve the optimal rate without mini-batches, we further develop a batch-free method that incorporates a first-order approximation and a clipping technique for function value estimation.
Finally, we validate the effectiveness of our proposed methods through experiments on risk-averse portfolio optimization and hierarchical tilted empirical risk minimization. 
\end{abstract}

\section{Introduction}
In this paper, we study stochastic multi-level compositional optimization~\cite{Yang2019MultilevelSG,chen2021solving,Zhang2021MultiLevelCS,doi:10.1137/21M1406222,TPAMI:2025:Jiang} of the form
\begin{align}
    \min_\x F(\x) = f_T \circ \cdots \circ f_{1} (\x).
\end{align} 
Here, the objective is a nested composition of $T$ smooth and non-convex functions.
We assume that only unbiased stochastic estimates are available for the function values and gradients at each level $i$, denoted as $ f_i(\cdot;\xi_i)$ and $\nabla f_i(\cdot;\xi_i)$, such that $\E_{\xi_i}\left[f_i(\cdot;\xi_i)\right] = f_i(\cdot)$ and $\E_{\xi_i}\left[\nabla f_i(\cdot;\xi_i)\right] = \nabla f_i(\cdot)$.

Multi-level optimization enjoys broad applications in machine learning, including risk-averse optimization~\cite{Bruno2016RiskNA,Shapiro2009LecturesOS}, reinforcement learning~\cite{Dann2014PolicyEW}, model-agnostic meta-learning~\cite{Ji2020MultiStepMM}, graph neural network training~\cite{10.1145}, stochastic variational inequalities~\cite{doi:10.1137/20M1312952}, adversarial attacks and defenses in deep neural networks~\cite{8611298}, risk management~\cite{Cole2013HowDR,Dentcheva2015StatisticalEO}, and hierarchical tilted empirical risk minimization~\cite{li2020tilted,li2021tilted}. 
Unlike traditional single-level stochastic optimization, where unbiased stochastic gradient evaluations are readily available, obtaining an unbiased estimate of the overall gradient and function value is generally infeasible in multi-level problems.
This difficulty stems from the nested composition structure, specifically:
\begin{align*}
    \E_{\xi_T,\cdots,\xi_2,\xi_1} [f_T( \cdots f_2(   f_1( \x;\xi_1);\xi_2)\cdots ;\xi_T)]
    \neq  \E_{\xi_T} [f_T(\cdots\E_{\xi_{2}} [f_{2}(  \E_{\xi_1} [f_1( \x;\xi_1)];\xi_{2})]\cdots ;\xi_T)].
\end{align*}
Furthermore, estimation errors can accumulate as the composition depth increases, making the problem more challenging.
To improve sample complexity, existing literature often imposes additional \emph{average smoothness} and \emph{average Lipschitz} conditions~\cite{chen2021solving,Zhang2021MultiLevelCS,TPAMI:2025:Jiang}, which require that stochastic estimates of each layer be smooth and Lipschitz continuous.
However, in single-level stochastic optimization, it is well known that momentum-based methods can achieve the optimal sample complexity of $\mathcal{O}(\epsilon^{-4})$ under standard smoothness assumptions (rather than the stronger average smoothness condition).
This motivates the following question: \emph{Is it possible to obtain the optimal rate in stochastic multi-level optimization without these additional assumptions?}

In this work, we answer this question affirmatively by proving the optimal $\mathcal{O}(\epsilon^{-4})$ sample complexity using a simple momentum-based approach with mini-batches.
Our analysis avoids the restrictive average smoothness and Lipschitz continuous conditions. Furthermore, we provide a normalized variant that eliminates the need for prior knowledge of problem-dependent constants when setting hyperparameters.
Finally, to avoid large batches, we develop a batch-free method that leverages first-order approximation together with clipping, preserving the optimal sample complexity with a constant batch size.
To summarize, our main contributions are as follows:
\begin{itemize}
\item Our proposed MML method uses momentum-based estimators to track function values and gradients, obtaining the optimal $\mathcal{O}(\epsilon^{-4})$ sample complexity with mini-batches. Our analysis does not rely on the average smoothness and Lipschitz continuity used in previous literature.
\item Employing a normalization technique, our MML-v2 removes the reliance on problem-dependent constants to set hyperparameters while preserving optimal sample complexity.
\item To avoid large batches, the proposed MMLFA method uses a first-order approximation for function-value tracking, combined with a clipping operation to control the estimator range. This design allows us to maintain the optimal rate with a constant batch size.
\item We validate the effectiveness of the proposed methods through numerical experiments on risk-averse portfolio optimization and hierarchical tilted empirical risk minimization.
\end{itemize}

\section{Related Work}
In this section, we review related work on stochastic multi-level optimization and its decentralized/federated variants.
\subsection{Multi-level Optimization}
Multi-level compositional optimization was first investigated by \citet{Yang2019MultilevelSG}, who developed stochastic methods based on an extrapolation-interpolation scheme. 
They established a sample complexity of $\mathcal{O} ( \epsilon^{-(7+T)/2})$ for non-convex objectives and $\mathcal{O} ( \epsilon^{-(3+T)/4} )$ for strongly convex objectives.
To further improve the convergence rate, \citet{chen2021solving} proposed the stochastically corrected stochastic compositional gradient method, deriving a sample complexity of $\mathcal{O}(\epsilon^{-4})$.
However, their analysis requires the stronger average smoothness assumption, which is typically used to obtain improved $\mathcal{O}(\epsilon^{-3})$ rates via variance reduction.
Concurrently, \citet{doi:10.1137/21M1406222} obtained $\mathcal{O}(\epsilon^{-6})$ sample complexity for a momentum-based method with mini-batches and improved the rate to $\mathcal{O}(\epsilon^{-4})$ using linearized stochastic estimates.
However, their analysis neglects dependency issues, rendering the results invalid (a detailed discussion is provided in the appendix).
To enhance the rate under the average smoothness assumption, \citet{Zhang2021MultiLevelCS} introduced the normalized proximal approximate gradient method, which evaluates gradients via nested variance reduction with large batch sizes.
They achieved a total sample complexity of $\mathcal{O}(\epsilon^{-3})$ in the expectation setting, and $\mathcal{O}(n+\sqrt{n}\epsilon^{-2})$ for the finite-sum structure, where $n$ denotes the number of component functions across all levels.
Similar results were obtained by \citet{jiang2022optimal, TPAMI:2025:Jiang} using a constant batch size, further deriving optimal rates for convex and strongly convex functions.

Beyond smooth non-convex settings, \citet{doi:10.1137/20M1312952} considered non-smooth non-convex multi-level optimization, assuming that functions at each level are locally Lipschitz continuous and differentiable in a generalized sense.
For smooth and convex problems, \citet{10.1007/s10107-024-02090-3} proposed a stochastic sequential dual method, achieving a sample complexity of $\mathcal{O}(\epsilon^{-2})$.
Constrained multi-level optimization has also been studied recently. 
\citet{NEURIPS2022_7e16384b} proposed a projection-free conditional gradient-type algorithm, attaining $\mathcal{O}(\epsilon^{-2})$ and $\mathcal{O}(\epsilon^{-3})$ complexities for stochastic first-order and linear-minimization oracles.
Following this, \citet{ICML:2024:Jiang} presented a comprehensive study on projection-free methods, covering non-convex, convex, and strongly convex objectives.
Recent studies have also focused on analyzing stability and generalization performance~\cite{pan2024stability} and extending multi-level optimization to minimax settings~\cite{zhang2025on}.

\subsection{Decentralized/Federated Multi-level Optimization}
Multi-level optimization has also been widely studied in decentralized and federated environments.
\citet{pmlr-v238-gao24b} first studied decentralized stochastic multi-level optimization, proposing two algorithms that achieve level-independent convergence rates for non-convex problems under milder conditions.
For smooth and strongly convex objectives, \citet{10507158} proposed the federated stochastic compositional gradient extrapolation method. 
Based on a double-loop strategy, this method attains an optimal $\mathcal{O}(\log (1/\epsilon ))$ communication complexity while maintaining an almost optimal $\widetilde{ \mathcal{O}}( 1/\epsilon )$ sample complexity.
To address imbalanced data classification in federated learning, \citet{pmlr-v235-zhang24aw} developed a federated multi-level compositional minimax algorithm for AUC maximization, equipped with rigorous theoretical guarantees.
To achieve linear speedup under heterogeneous settings, \citet{pmlr-v235-gao24a} developed a federated stochastic multi-level algorithm based on a Jacobian-vector product estimator and proved level-independent convergence for non-convex objectives. 

\section{Momentum-based Multi-level Algorithm}
In this section, we first outline the assumptions used and then present the proposed momentum-based methods, followed by their convergence guarantees.
\subsection{Assumptions}
We adopt the following assumptions, which are standard in stochastic multi-level optimization~\cite{Yang2019MultilevelSG,chen2021solving,Zhang2021MultiLevelCS,doi:10.1137/21M1406222,TPAMI:2025:Jiang}.
\begin{ass} \label{asm:0} 
Denote $F_\star=\min_{\x} F(\x)$ and suppose $F\left(\x^{1}\right)- F_\star \leq \Delta_{F}$ for the initial point $\x^{1}$.
\end{ass}
\begin{ass} \label{asm:1} 
All functions $f_1, \dotsc, f_T$ are $L_f$-Lipschitz continuous, and their gradients $\nabla f_1, \dotsc, \nabla f_T$ are also $L_J$-Lipschitz continuous.
\end{ass}
\textbf{Remark:} Under Assumption~\ref{asm:1}, the objective $F$ is $L_F$-smooth, where
$L_F \coloneqq L_f^{2T-1}L_J\sum_{i=1}^T{L_f^{-i}}$.
\begin{ass}\label{asm:stochastic2} 
The function and gradient estimates are unbiased with bounded variance, i.e.,
\begin{gather*}
       \E_{\xi}\left[f_i(\x;\xi)\right] = f_i(\x), \quad
	  \E_{\xi}\left[\nabla f_i(\x;\xi)\right] = \nabla f_i(\x),\\
\E_{\xi}\left[\Norm{f_i(\x;\xi) - f_i(\x)}^2 \right]\leq \sigma^2, \quad
	\E_{\xi}\left[\Norm{\nabla f_i(\x;\xi) - \nabla f_i(\x) }^2 \right] \leq  \sigma^2.
\end{gather*}
\end{ass}
\textbf{Remark:} These assumptions are standard in traditional stochastic compositional optimization. Some existing multi-level literature~\cite{chen2021solving,Zhang2021MultiLevelCS,TPAMI:2025:Jiang} requires the stronger average smoothness and Lipschitz condition, namely, for each $1\le i\le T$:
\begin{align*}
	\E_{\xi}\left[\Norm{f_i(\x;\xi)-  f_i(\y;\xi)}^2\right] &\leq \LL_f^2 \Norm{\x - \y}^2,\,
		\E_{\xi}\left[\Norm{\nabla f_i(\x;\xi)- \nabla f_i(\y;\xi)}^2\right]\leq  \LL_J^2 \Norm{\x - \y}^2.
\end{align*}
However, this assumption is not required in our work.

\subsection{The Proposed Method}
We now present our algorithm for stochastic multi-level optimization. 
In traditional stochastic non-convex optimization~\cite{kingma:adam,cutkosky2019momentum}, we can directly access an unbiased stochastic gradient estimate of the overall objective $\nabla F(\cdot)$.
However, this is no longer possible in multi-level optimization, even though unbiased estimates of the function and gradient at each individual layer are available. 
Since the overall gradient $\nabla F(\x)$ typically involves nested expectations, multiplying the stochastic gradients of each layer via the chain rule, i.e., $\prod_{i=1}^T \nabla f_i (\cdot;\xi_i)$, yields only a biased estimation, potentially resulting in significant estimation error.
This limitation arises from the nested nature of the problem, since $\E_{\xi,\zeta} [f_2(f_1(\x;\zeta);\xi)] \neq \E_{\xi}[f_2(\E_{\zeta}[f_1(\x;\zeta)]);\xi]$.

To address this issue, we first provide a detailed analysis of the gradient estimation error. For clarity, we consider the simple two-level case where $T=2$.
The objective function is $F(\x) = \E_\xi[f_2( \E_\zeta[f_1(\x;\zeta)]; \xi  ) ]$, and its gradient is $\nabla F(\x) = \E_\xi[\nabla f_2( \E_\zeta[f_1(\x;\zeta)]; \xi  ) ]\cdot  \E_\omega[\nabla f_1(\x;\omega)]$.
Note that there exist three distinct expectations $\E_\xi,\E_\zeta, \E_\omega$. 
To construct an effective gradient estimator, we must derive accurate estimates of the gradients $\E_\xi[\nabla f_2(\cdot;\xi) ]\cdot\E_\omega[\nabla f_1(\cdot;\omega)]$ and the inner function values $ \E_\zeta[f_1(\cdot;\zeta)]$.
Therefore, we first employ a momentum-based estimator to track the function values of each  layer $i$ at time step $k$:
\begin{align}\label{u}
    \u_i^k = (1-\alpha)\u_{i}^{k-1}  + \alpha \frac{1}{B_1}\sum_{j=1}^{B_1}f_i(\u^k_{i-1};\xi^k_{i,j}),
\end{align}
where $\u^k_{0}=\x^k$ and $\u^1_i = \frac{1}{B_0}\sum_{j=1}^{B_0}f_i(\u^1_{i-1};\xi^1_{i,j})$.
Note that multi-level problems are more difficult than two-level problems, because the estimation errors may accumulate exponentially with the composition depth. This is why we employ the level-wise function estimators $\u_i^k$. We analyze the coupled recursions for $\Gamma_i$ and $\Lambda_i$ in the following lemma, where we denote the error of the function estimator as $\Gamma_i = \frac{1}{K}\sum_{k=1}^K  \E\left[\Norm{\u^k_i -  f_i(\u^k_{i-1})}^2\right]$ and the difference between adjacent steps as $\Lambda_i=\frac{1}{K}\sum_{k=1}^K \E\left[\Norm{\u_{i}^{k+1} - \u_{i}^{k}}^2\right]$. We will then use them to control the gradient estimator error. 

\begin{algorithm}[!tb]
	\caption{MML method}
	\label{alg:1}
	\begin{algorithmic}[1]
 \STATE {\bfseries Input:} parameters $\eta$, $\alpha$ 
		\FOR{time step $k = 1$ {\bfseries to} $K$}
		\STATE Set $\u^k_0 = \x^k$
\STATE Compute the function estimator according to equation~(\ref{u})
        \STATE Compute the gradient estimator $\z^k$ according to equation~(\ref{z})
\STATE (v1) Update the weight: $\x^{k+1} = \x^k - \eta \z^k $
		\STATE (v2) Update the weight: $\x^{k+1} = \x^k - \eta \frac{\z^k}{\Norm{\z^k}} $
		\ENDFOR
	\end{algorithmic}
\end{algorithm}
\begin{lemma}\label{L1} 
Set $B_1 \leq \alpha^2 B_0 K$. We can ensure that
\begin{align*}
\Gamma_i \leq \frac{2\LL_f^2}{\alpha^2 }  \Lambda_{i-1} + \frac{2\alpha   \sigma^2}{B_1};\quad
 \Lambda_i \leq {2\alpha^2\LL_f^2}\Lambda_{i-1} + \frac{\alpha^2\sigma^2}{B_1}
 + 2\alpha^2 \Gamma_i.
\end{align*}
\end{lemma} 
The function estimation error $\Gamma_i$ and the step difference $\Lambda_i$ are mutually dependent. We can bound them simultaneously in the following lemma.
\begin{lemma}\label{L2} 
Let $C_{i+1} = 5+6\LL_f^2 C_i$ with $C_0=1$. Then
\begin{align*}
\Gamma_i  \leq &C_i \left(\frac{1}{\alpha^2 K}  \sum_{k=1}^K \E\left[\Norm{\x^{k+1} - \x^{k}}^2\right]  + \frac{ \sigma^2}{B_1} \right);\\
\Lambda_i \leq & C_i  \left(\frac{1}{K}  \sum_{k=1}^K \E\left[\Norm{\x^{k+1} - \x^{k}}^2\right]  + \frac{\alpha^2 \sigma^2}{B_1} \right).
\end{align*}
\end{lemma}	
\textbf{Remark:} 
By properly setting the learning rate, we can control the term $\E\left[\Norm{\x^{k+1} - \x^{k}}^2\right]$. Consequently, using a sufficiently large batch size $B_1$ yields a small function-value estimation error.

Using the function value estimators $\u_i^k$, we develop the overall gradient estimator via the chain rule:
\begin{align}\label{z}
    \z^k = & (1-\alpha)\z^{k-1}  +\alpha \frac{1}{B_1} \sum_{j=1}^{B_1}  \prod_{i=1}^T \nabla f_i(\u^{k}_{i-1};\xi^{k}_{i,j}), 
\end{align}
where $\z^1 = \frac{1}{B_0} \sum_{j=1}^{B_0}  \prod_{i=1}^T \nabla f_i(\u^{1}_{i-1};\xi^{1}_{i,j})$. 
This gradient estimator also follows a momentum update, and its error is bounded as follows.
\begin{lemma}\label{L3}
We have the following guarantee:
    \begin{align*}
    \frac{1}{K}\sum_{k=1}^K \E\left[\Norm{\z^{k}-\nabla F(\x^{k})}^2\right]
    \leq & \frac{2\alpha }{B_1} T^2(\sigma^2 +L_f^2 )^{T}+\frac{ \sigma^2}{B_1} 4T L_F^2\sum_{i=1}^{T-1} C_i  \\
    & +\left(1+T\sum_{i=1}^{T-1} C_i \right) \frac{4L_F^2}{\alpha^2 K}\sum_{k=1}^{K+1} \E\left[\Norm{\x^{k+1} - \x^{k}}^2\right].
\end{align*}
\end{lemma}
\textbf{Remark:} With suitable choices of the learning rate $\eta$ and momentum parameter $\alpha$, we can ensure a small gradient estimation error using a large batch size $B_1$.

Once the gradient estimator is obtained, the decision variable is updated using either standard SGD (i.e., $\x^{k+1} = \x^k -\eta \z^k$) or normalized SGD (i.e., $\x^{k+1} = \x^k -\eta \frac{\z^k}{\Norm{\z^k}}$). The complete procedure, named the Momentum-based Multi-Level method (MML), is presented in Algorithm~\ref{alg:1}.
We now present the convergence guarantees below.
\begin{theorem}\label{T1}
Under Assumptions~\ref{asm:0}, \ref{asm:1}, \ref{asm:stochastic2}, by setting $B_0 = \mathcal{O}(1)$, $B_1 = \mathcal{O}(K)$, $\eta = \mathcal{O}(1)$ and $\alpha =\mathcal{O}(1)$, our MML-v1 algorithm ensures 
\begin{align*}
    \frac{1}{K} \sum_{k=1}^K \E\left[\Norm{\nabla F(\x^k)}\right] \leq \mathcal{O}\left(\frac{1}{\sqrt{K}}\right).
\end{align*}
\end{theorem}
\textbf{Remark:} To ensure that $\frac{1}{K} \sum_{k=1}^K \E\left[\Norm{\nabla F(\x^k)}\right] \leq \epsilon$, the iteration complexity is $K=\mathcal{O}(\epsilon^{-2})$ and the overall sample complexity is $B_1 K=\mathcal{O}(\epsilon^{-4})$. 
This matches the lower bound for stochastic single-level optimization under standard smoothness assumptions~\citep{Arjevani2019LowerBF} and is therefore optimal. Notably, this is the first result to obtain optimal $\mathcal{O}(\epsilon^{-4})$ sample complexity for multi-level optimization under standard smoothness assumptions.

While Theorem \ref{T1} attains optimal sample complexity, it requires knowledge of problem-dependent constants to set hyperparameters, e.g., $\eta \leq \min\{\frac{1}{2L_F}, \frac{\alpha}{(4L_F\sqrt{1+T\sum_{i=1}^{T-1} C_i})}\}$. 
We next show that this requirement can be removed using the normalization technique, as shown in the following theorem.
\begin{theorem}\label{T2}
Under Assumptions~\ref{asm:0}, \ref{asm:1}, \ref{asm:stochastic2}, by setting $B_0=1$, $B_1 = \lfloor K/4 \rfloor$, $\eta = K^{-1/2}$, and $\alpha = 1/2$, our MML-v2 algorithm ensures 
\begin{align*}
    \frac{1}{K} \sum_{k=1}^K \E\left[\Norm{\nabla F(\x^k)}\right] \leq \mathcal{O}\left(\frac{1}{\sqrt{K}}\right).
\end{align*}
\end{theorem}
\textbf{Remark:} Theorem~\ref{T2} provides the same rate while avoiding the use of problem-dependent constants for setting hyperparameters. Note that $\alpha = 1/2$ is not essential to the analysis and can be any positive constant in the range $(0,1)$. A similar normalization technique is also used by \citet{pmlr-v238-hubler24a} in single-level problems. However, our analysis is different from \citet{pmlr-v238-hubler24a}, which analyzes with decreasing $\eta_t$ and increasing $\alpha_t$, and obtains a nearly optimal rate of $\mathcal{O}(K^{-1/4}\log K)$.

\section{Momentum-based Multi-level Optimization with First-order Approximation}
In the previous section, we achieved optimal sample complexity using mini-batches.
While large batches can be beneficial (e.g., enabling parallelism), they would become impractical when hardware or memory constraints prevent the use of $B_1=\mathcal{O}(K)$. 
In this section, we introduce an alternative algorithm that achieves the optimal sample complexity without relying on large batches.

A careful examination of the analysis reveals that the large batch requirement stems from Lemma~\ref{L2}, specifically the term ${\sigma^2}/{B_1}$ in the estimation error.
To ensure a small function estimation error, a large $B_1$ is required, rather than a constant batch size of $B_1 = \mathcal{O}(1)$.
To address this issue, we propose to track function values using the first-order Taylor approximation:
\begin{align*}
    f_i(\u_{i-1}^k)\approx f_i(\u_{i-1}^{k-1})+\langle \nabla f_i(\u_{i-1}^{k}), \u_{i-1}^k-\u_{i-1}^{k-1}\rangle.
\end{align*}
Recall that the estimator $\u_i^k$ tracks $f_i(\u_{i-1}^k)$ and $\u_i^{k-1}$ estimates $f_i(\u_{i-1}^{k-1})$.
Moreover, we can obtain a stochastic estimate of $\nabla f_i(\u_{i-1}^{k})$, denoted by $\nabla f_i(\u_{i-1}^{k};\xi_{i}^k)$. Consequently, we approximate
\begin{align*}
    \u_i^k \approx \u_i^{k-1}+\langle \nabla f_i(\u_{i-1}^{k};\xi_{i}^k), \u_{i-1}^k-\u_{i-1}^{k-1}\rangle.
\end{align*}
Since $f_i(\u_{i-1}^k;\xi_i^k)$ is an unbiased estimator of $f_i(\u_{i-1}^k)$, we combine these two estimators to form
\begin{equation}\label{u+}
    \begin{split}
        \tilde{\u}_i^k  =\alpha f_i(\u_{i-1}^k;\xi_i^k)
  + (1-\alpha)\left( \u_i^{k-1}+\langle \nabla f_i(\u_{i-1}^{k};\xi_{i}^k), \u_{i-1}^k-\u_{i-1}^{k-1}\rangle \right),
    \end{split}
\end{equation}
where $\tilde{\u}^1_i = f_i(\u^1_{i-1};\xi^1_{i})$.
The parameter $\alpha$ serves as a trade-off factor between two estimation methods.
Although this new estimator requires computing $\nabla f_i(\u_{i-1}^{k};\xi_{i}^k)$, this term is also needed for constructing the gradient estimator. 
Thus, we can improve the estimator without increasing the overall sample complexity.
Compared to the traditional momentum-based estimator, we incorporate the additional term $(1-\alpha)\langle \nabla f_i(\u_{i-1}^{k};\xi_{i}^k), \u_{i-1}^k-\u_{i-1}^{k-1}\rangle$, which compensates for the drift $\u_{i}^k-\u_{i}^{k-1}$ via the first-order Taylor expansion.
After computing $\tilde{\u}_i^k$, we apply a clipping operation to keep the estimator within a controlled range:
\begin{align*}
    \u_i^k = \operatorname{Clip}_{C}[\tilde{\u}_i^k].
\end{align*}
This operation is straightforward to implement:
\begin{align*}
    \u_i^k=
    \begin{cases}
        \tilde{\u}_i^k, & \text{if } \Norm{\tilde{\u}_i^k}\le C,\\[3pt]
        \displaystyle \frac{C}{\Norm{\tilde{\u}_i^k}}\tilde{\u}_i^k, & \text{otherwise}.
    \end{cases}
\end{align*}
Here, parameter $C$ is the clipping threshold and we require $\Norm{f_i(\u_{i-1}^k)} \leq C$ for each $i$ and $k$. In practice, we simply set $C$ as a large value, as discussed in our experiments.
Next, we derive the guarantee for this new estimator.
As before, we define that $\Gamma_i = \frac{1}{K}\sum_{k=1}^K  \E\left[\Norm{\u^k_i -  f_i(\u^k_{i-1})}^2\right]$ and  $\Lambda_i=\frac{1}{K}\sum_{k=1}^K \E\left[\Norm{\u_{i}^{k+1} - \u_{i}^{k}}^2\right]$.
\begin{lemma}\label{L1+} 
Set $1 \leq \alpha^2  K$, and we can ensure that
\begin{align*}
	    \Gamma_i \leq &\frac{(4L_f^2+2L_JC +2\sigma^2)}{\alpha}  \Lambda_{i-1} + 3\alpha   \sigma^2;\quad
 \Lambda_i \leq  2\alpha^2\sigma^2  + (6\LL_f^2+2\sigma^2)\Lambda_{i-1} +  (4\alpha^2)\Gamma_i.
\end{align*}
\end{lemma}	
\textbf{Remark:} 
A key improvement over  Lemma~\ref{L1} is the reduction of the factor for $\Lambda_{i-1}$ from ${1}/{\alpha^2}$ to ${1}/{\alpha}$ in the recurrence of $\Gamma_i$. 
This allows us to use a much smaller $\alpha$ to control the term $\alpha \sigma^2$, thereby removing the reliance on a large batch size $B_1$.
We highlight this difference in the following lemma.
\begin{algorithm}[t]
	\caption{MMLFA method}
	\label{alg:2}
	\begin{algorithmic}[1]
	\STATE {\bfseries Input:} parameters $\eta$, $\alpha$, $C$ 
		\FOR{time step $k = 1$ {\bfseries to} $K$}
		\STATE Set $\u^k_0 = \x^k$
  \STATE Compute the function estimator $\tilde{\u}_i^k$ according to equation~(\ref{u+})
  \STATE Conduct the clipping operation $\u_i^k = \operatorname{Clip}_{C}[\tilde{\u}_i^k]$
    \STATE Compute the gradient estimator $\z^k$ according to equation~(\ref{z+})
\STATE Update the weight: $\x^{k+1} = \x^k - \eta \z^k $
		\ENDFOR
	\end{algorithmic}
\end{algorithm}
\begin{lemma}\label{L2+} 
Let $D = 2(7L_f^2+2L_JC +3\sigma^2)$. We obtain:
\begin{align*}
\Gamma_i  \leq & D^i \frac{1}{\alpha K}\sum_{k=1}^K\E\left[\Norm{\x^{k+1} - \x^{k}}^2\right]  + (8 D^{i} +3)\alpha   \sigma^2;\\
\Lambda_i  \leq & D^i \frac{1}{K}\sum_{k=1}^K\E\left[\Norm{\x^{k+1} - \x^{k}}^2\right] + 8\alpha^2   \sigma^2 D^i.
\end{align*}
\end{lemma}	
\textbf{Remark:} Compared to Lemma~\ref{L2}, the second term in the function estimation error changes from $\mathcal{O}(\sigma^2/B_1)$ to $\mathcal{O}(\alpha \sigma^2)$, which can be controlled by a small $\alpha$ rather than a large batch size. Also, the factor of the first term is improved from $1/\alpha^2$ to $1/\alpha$, enabling the use of much smaller $\alpha$.

We then construct the overall gradient estimator using these $\u_i^k$ values without batches:
\begin{align}\label{z+}
    \z^k = & (1-\alpha)\z^{k-1}  +\alpha  \prod_{i=1}^T \nabla f_i(\u^{k}_{i-1};\xi^{k}_{i}), 
\end{align}
where $\z^1 =  \prod_{i=1}^T \nabla f_i(\u^{1}_{i-1};\xi^{1}_{i})$. Similarly, the gradient  estimation error can be bounded as follows.
\begin{lemma}\label{L3+}
Let $D_0 = 2 T^2{(\sigma^2+L_f^2)^T} +4 T L_F^2   \sigma^2(8 D^{T} +3T)$. We have the following guarantee:
    \begin{align*}
    \frac{1}{K}\sum_{k=1}^K \E\left[\Norm{\z^{k}-\nabla F(\x^{k})}^2\right] \leq D_0 \alpha
     + 4L_F^2\left(1+ T  D^T\right) \frac{1}{\alpha^2K}\sum_{k=1}^{K+1}\E\left[\Norm{\x^{k+1} - \x^{k}}^2\right].
\end{align*}
\end{lemma}               
\textbf{Remark:} 
By properly setting the learning rate $\eta$ and momentum parameter $\alpha$, we ensure a small gradient estimation error with a constant batch size, as $\frac{ \sigma^2}{B_1} 4T L_F^2\sum_{i=1}^{T-1} C_i$ in Lemma~\ref{L3} is eliminated.

Finally, we update the decision variable $\x$ using an SGD-style step. The complete procedure, named Momentum-based Multi-Level method with First-order Approximation (MMLFA), is described in Algorithm~\ref{alg:2}.
We now provide the convergence guarantee for MMLFA.
\begin{theorem}\label{T1+}
Under Assumptions~\ref{asm:0}, \ref{asm:1}, \ref{asm:stochastic2} and supposing that $\Norm{f_i(\u_{i-1}^k)} \leq C$, by setting $\eta = \mathcal{O}(K^{-1/2})$ and $\alpha =\mathcal{O}(K^{-1/2})$, our MMLFA algorithm ensures that:
\begin{align*}
    \frac{1}{K} \sum_{k=1}^K \E\left[\Norm{\nabla F(\x^k)}\right] \leq \mathcal{O}\left(\frac{1}{K^{1/4}}\right).
\end{align*}
\end{theorem}
\textbf{Remark:} 
To ensure $\frac{1}{K} \sum_{k=1}^K \E\left[\Norm{\nabla F(\x^k)}\right] \leq \epsilon$, the sample complexity is $K=\mathcal{O}(\epsilon^{-4})$. This rate matches the previous result while eliminating the need for large batches.

\section{Experiments}
In this section, we validate the effectiveness of the proposed methods on mean-deviation risk-averse portfolio optimization and hierarchical tilted empirical risk minimization.
We compare our methods with existing stochastic multi-level algorithms, including A-TSCGD~\cite{Yang2019MultilevelSG}, NLASG~\cite{doi:10.1137/21M1406222}, Nested-SPIDER~\cite{Zhang2021MultiLevelCS}, and SCSC~\cite{chen2021solving}.
Hyperparameters for each method are selected via grid search: the momentum parameter is tuned over the set $\{0.01, 0.05, 0.1, 0.5, 0.9, 0.95, 0.99\}$ and the learning rate is searched from the set $\{1, 10^{-1}, 10^{-2}, 10^{-3}, 10^{-4}\}$. 
We report the best results obtained from the grid search for each method. If the optimal hyperparameter lies on the boundary of the grid, we extend the search range to ensure it lies within the interior.
We set the batch size to 128 for each method.
For the clipping threshold $C$, we observe that a sufficiently large value works well in practice, avoiding the need for tuning. All results are averaged over 20 runs.

\subsection{Mean-deviation Risk-averse Portfolio Management}
\begin{figure*}[!t]
\begin{center}
\includegraphics[width=\rat\textwidth]{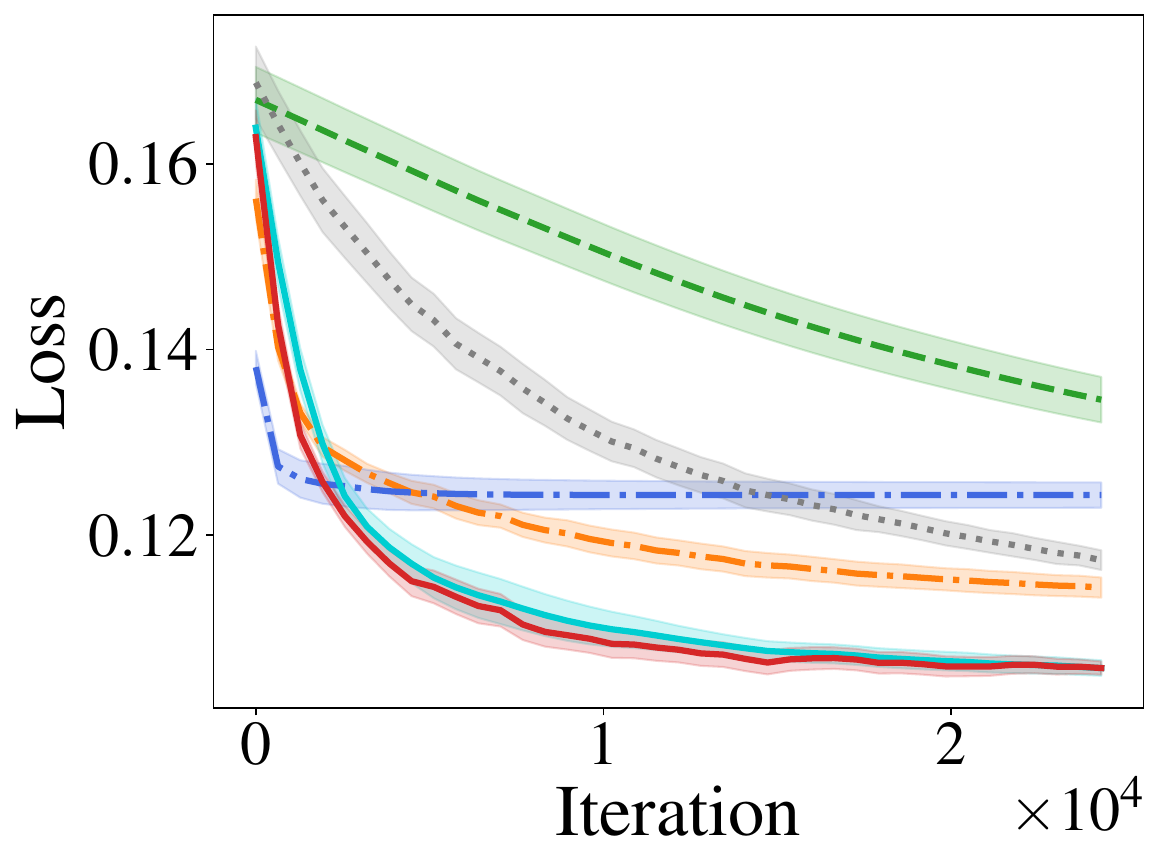}
\includegraphics[width=\rat\textwidth]{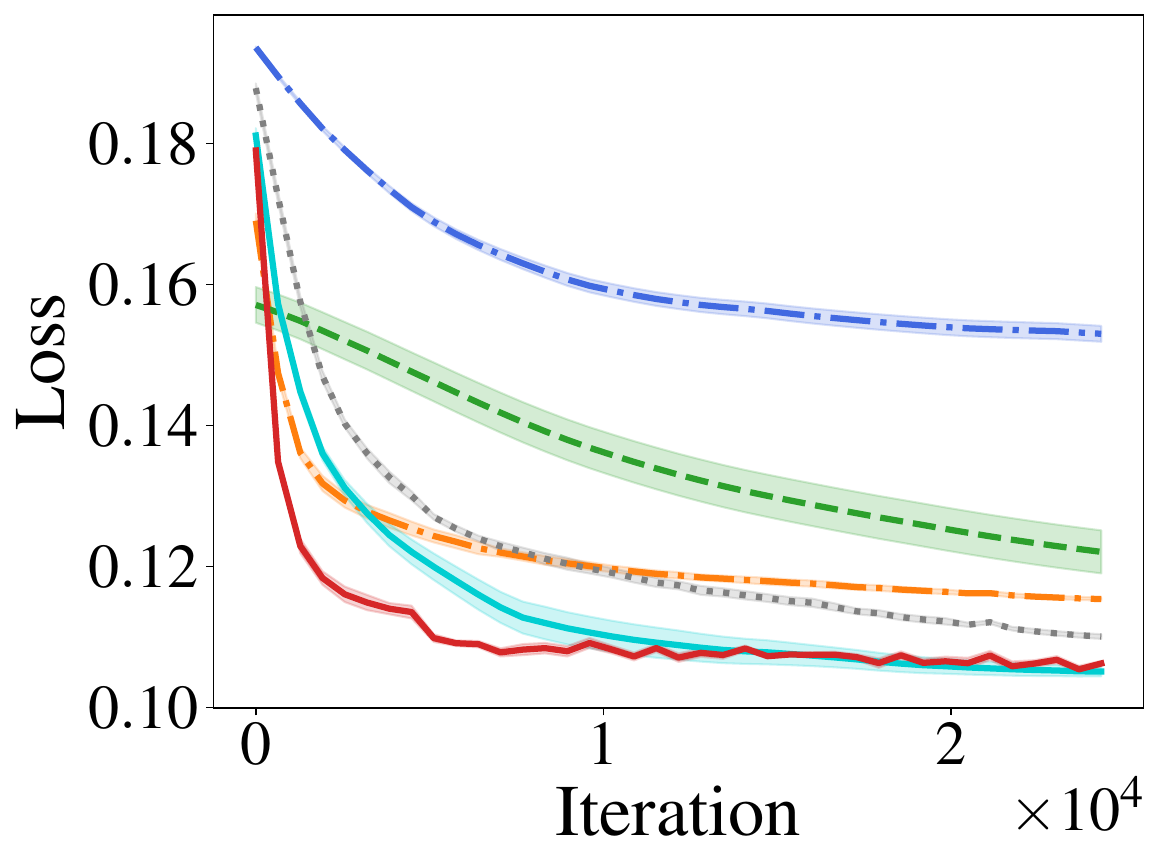}
\includegraphics[width=\rat \textwidth]{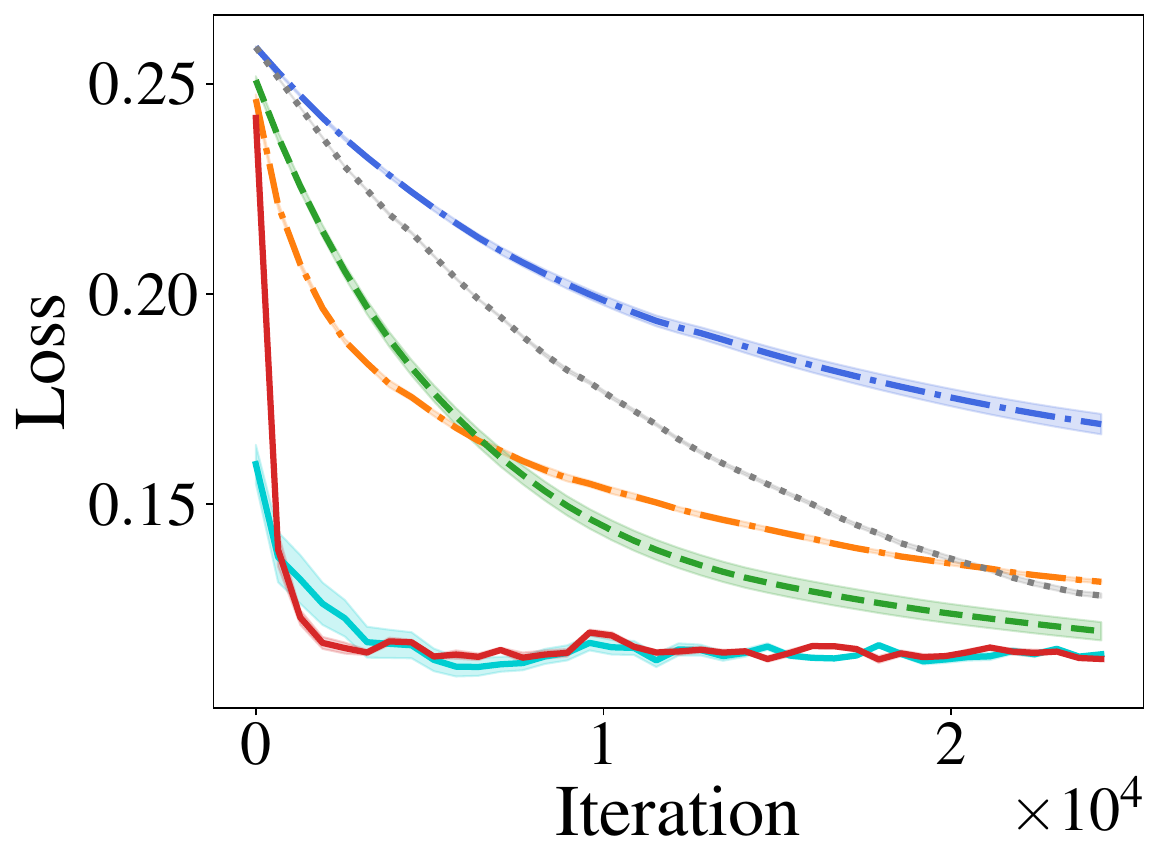}\\
\subfigure[Industry-10]{\includegraphics[width=\rat\textwidth]{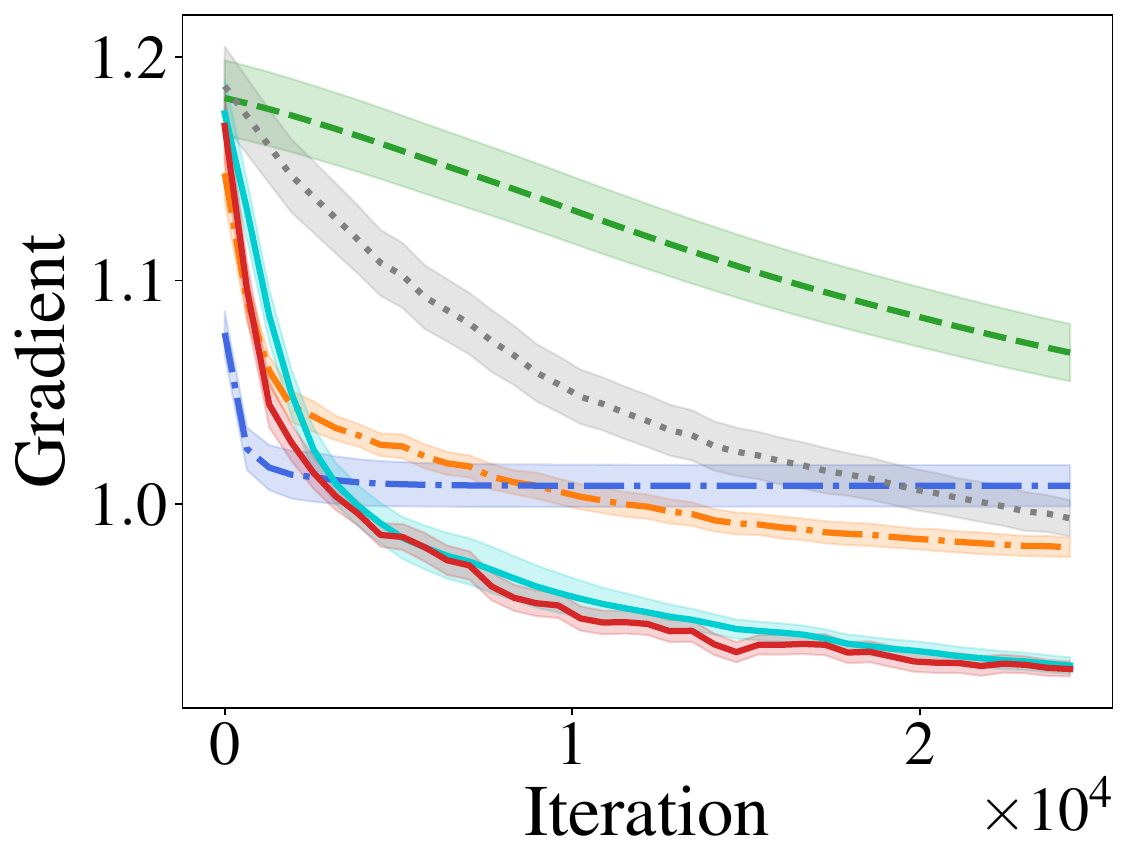}}
\subfigure[Industry-12]{\includegraphics[width=\rat\textwidth]{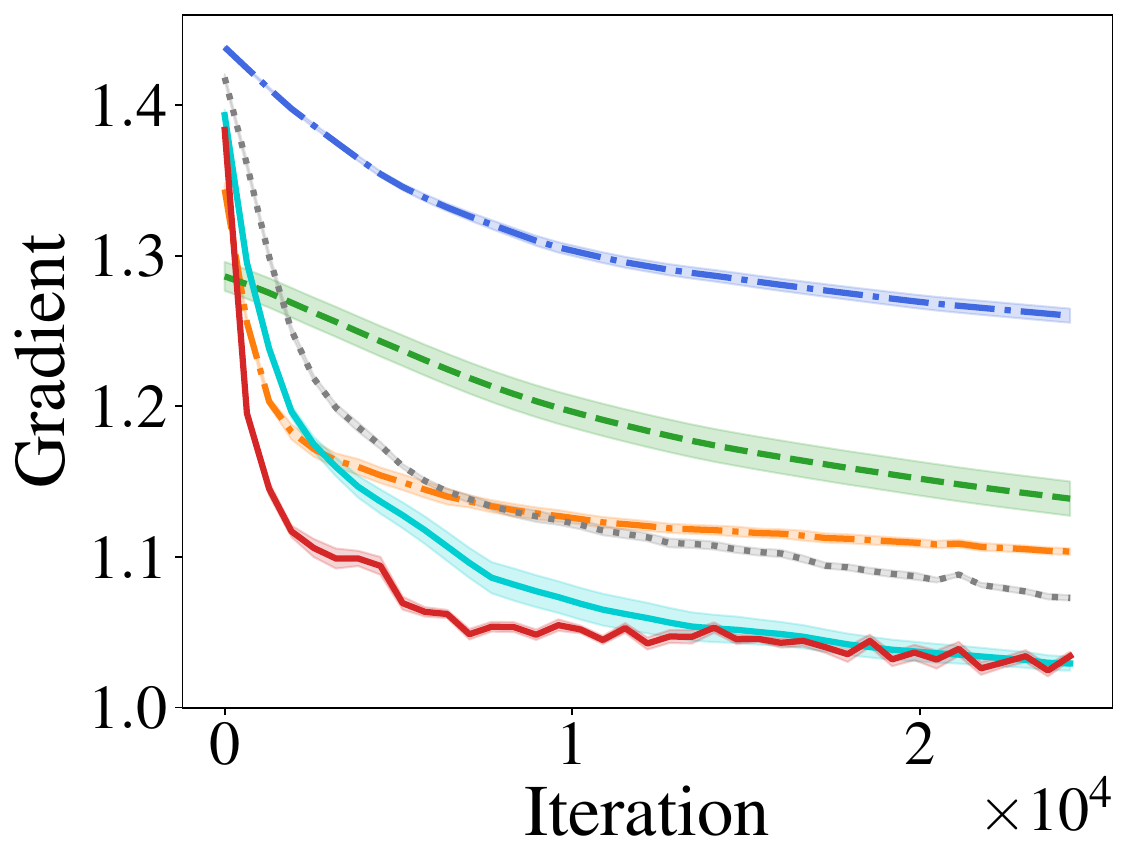}}
\subfigure[Industry-17]{\includegraphics[width=\rat\textwidth]{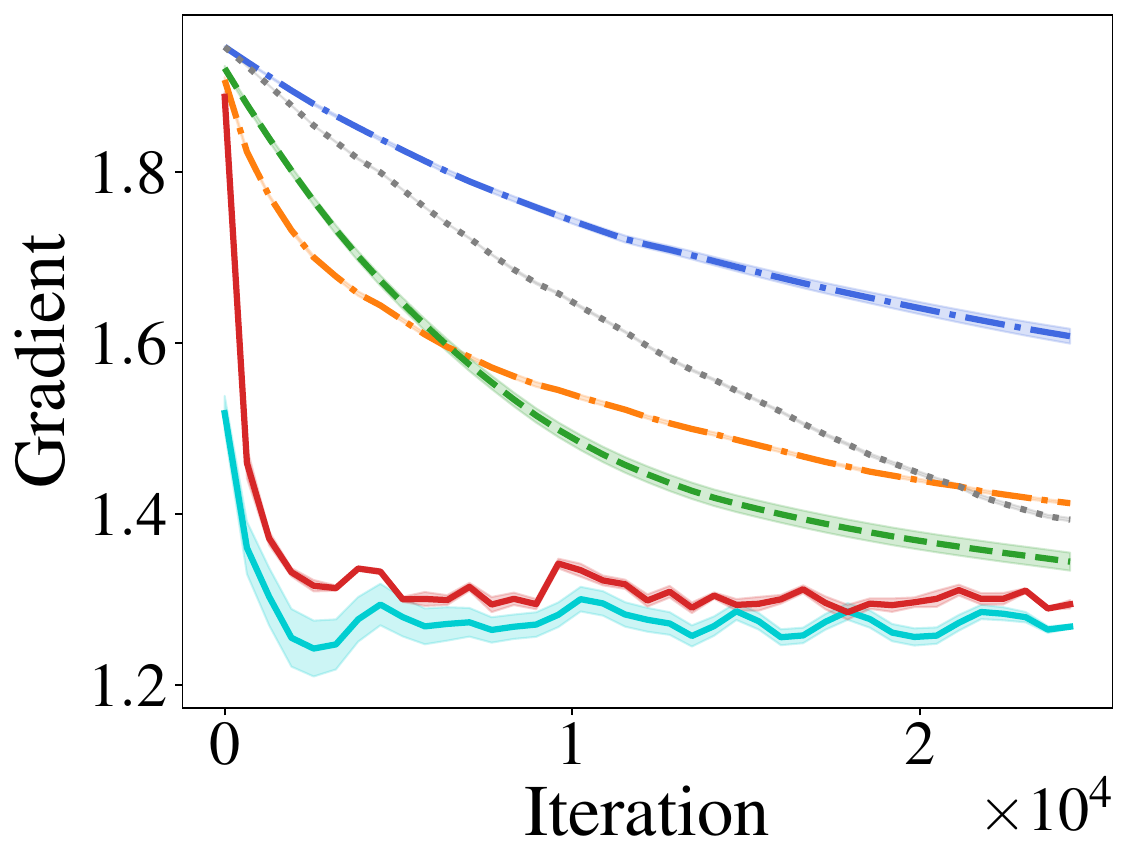}}
\includegraphics[width=0.9\textwidth]{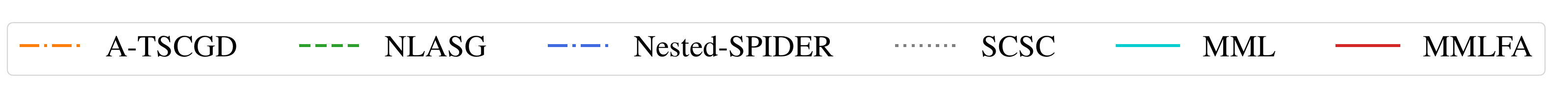}
\caption{Results for Risk-Averse Portfolio Optimization.}
		\label{fig:1}
	\end{center}
\end{figure*}

First, we investigate the risk-averse portfolio optimization.
Consider a portfolio of $d$ assets, where the vector $\r_{n} \in \mathbb{R}^{d}$ denotes the payoff of all assets at step $n \in \{1, \ldots, N\}$.
We aim to identify a strategy $\x$, where $[\x]_i$ is the fraction of total capital allocated to asset $i$. 
The goal is to maximize expected return while controlling risk. A robust formulation for this problem is the mean-deviation risk-averse optimization model~\citep{Shapiro2009LecturesOS}, in which risk is defined as the standard deviation.
This model is widely applied in practice and frequently serves as a benchmark in multi-level optimization research~\citep{Yang2019MultilevelSG, Zhang2020OptimalAF}. 
The problem is formulated as:
\begin{align*}
\min _{\x} -\frac{1}{N} \sum_{n=1}^{N}\left\langle \r_{n}, \x\right\rangle+\lambda \sqrt{\frac{1}{N} \sum_{n=1}^{N}\left(\left\langle \r_{n}, \x\right\rangle-\left\langle \bar{\r}, \x\right\rangle\right)^{2}},
\end{align*}
where $\bar{\r}= \frac{1}{N}\sum_{n=1}^{N} \r_{n}$. The above formulation constitutes a three-level stochastic compositional optimization problem, with levels defined as:
\begin{align*}
     f_{1}(\x)&=\left(\frac{1}{N} \sum_{n=1}^{N}\langle \r_{n}, \x\rangle, \x\right), \\
     f_{2}((z, \x))&=\left(z, \frac{1}{N} \sum_{n=1}^{N}\left(\left\langle \r_{n}, \x\right\rangle-z\right)^{2}\right),  \\
     f_{3}\left((y_{1}, y_{2})\right)&=-y_{1}+\lambda \sqrt{y_{2}}.
\end{align*}
In the experiments, we evaluate on three real-world datasets (Industry-10, Industry-12, and Industry-17) from the Kenneth R. French Data Library\footnote{https://mba.tuck.dartmouth.edu/pages/faculty/ken.french/}, which contain 10, 12, and 17 industrial assets, respectively.
Following \citet{Zhang2021MultiLevelCS}, we set $\lambda=0.2$.
Figure~\ref{fig:1} illustrates the loss values and gradient norms with respect to iterations~(since the batch size is fixed to 128, the sample complexity is just proportional to the iteration number).
Across all datasets, our methods (MML and MMLFA) demonstrate superior performance.
Both the loss and gradient norms decrease more rapidly and stabilize at lower levels compared to the other methods.
MML and MMLFA perform similarly overall: MMLFA performs slightly better on Industry-10/12, while MML is slightly better on Industry-17.

Note that although MML requires large batch sizes in the theoretical analysis, we do not observe significant performance degradation in the experiments when using smaller batches. In Appendix~\ref{addExp}, we validate the performance with different batch sizes, i.e., $\{48,64,96,128\}$. As shown in Figure~\ref{fig:01}, we do not observe significant performance degradation with smaller batches, demonstrating the practical robustness of MML with small batches. Overall, we would recommend using MML first, since it is simpler and therefore easier and more stable to implement. When the batch size is extremely small, MMLFA may be preferable. 

We further study the effect of the clipping threshold $C$ in MMLFA.
Theoretically, choosing $C$ above the true upper bound of the relevant function values does not change the order of the convergence rate.
Here, we vary $C$ within the range $\{0.5, 1, 5, 10, 10^{20}\}$, with results depicted in Figure~\ref{fig:2} in Appendix~\ref{addExp}. The results indicate that performance is stable provided $C$ is not too small (e.g., $C\ge 1$). 
In contrast, overly aggressive clipping (e.g., $C=0.5$) degrades performance, suggesting that a large default value of $C$ is reasonable in practice.

\subsection{Hierarchical Tilted Empirical Risk Minimization}
\begin{figure*}[!t]
	\begin{center}
\includegraphics[width=\rat\textwidth]{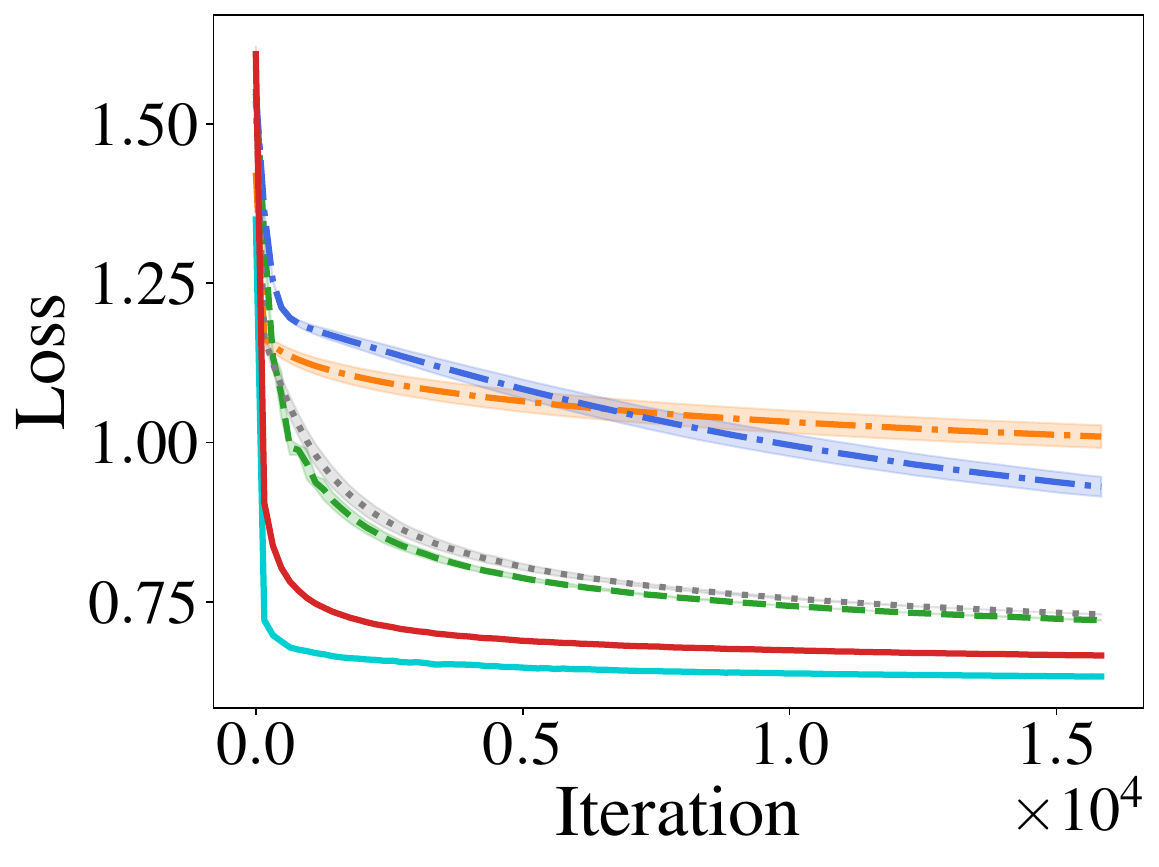}
\includegraphics[width=\rat\textwidth]{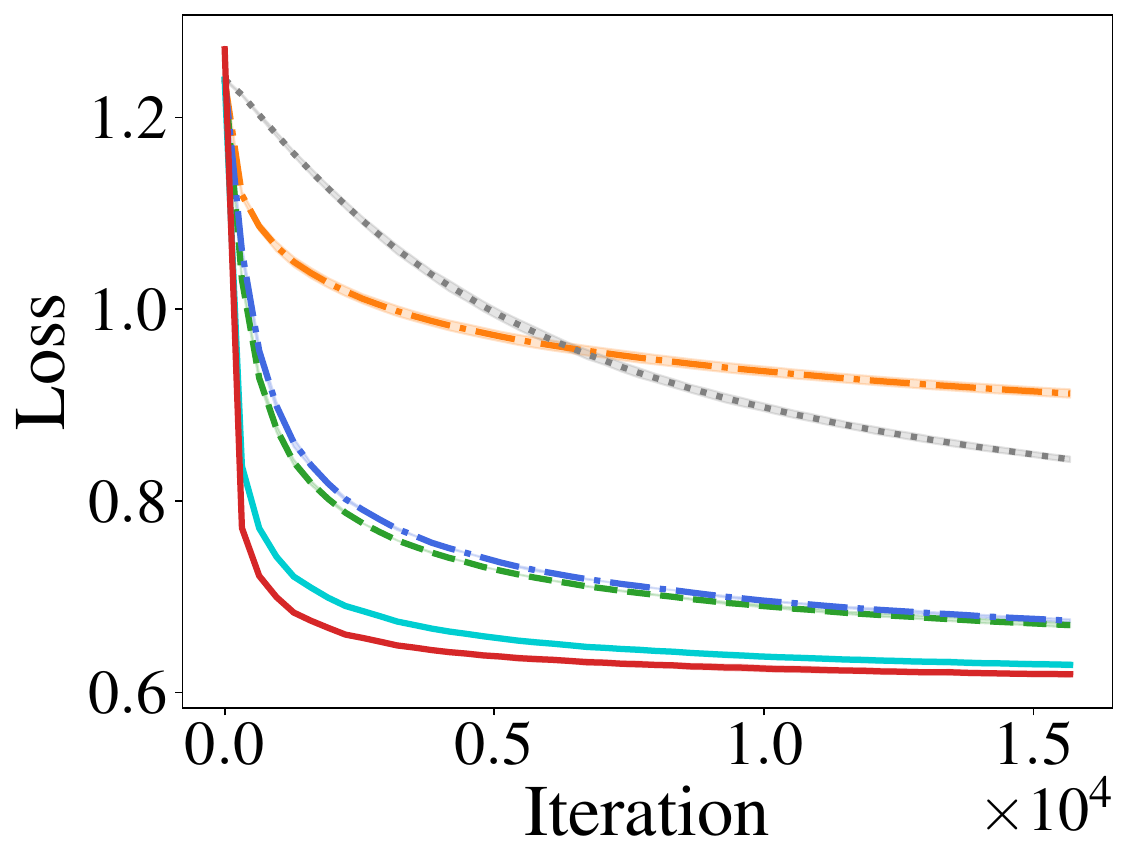}
\includegraphics[width=\rat\textwidth]{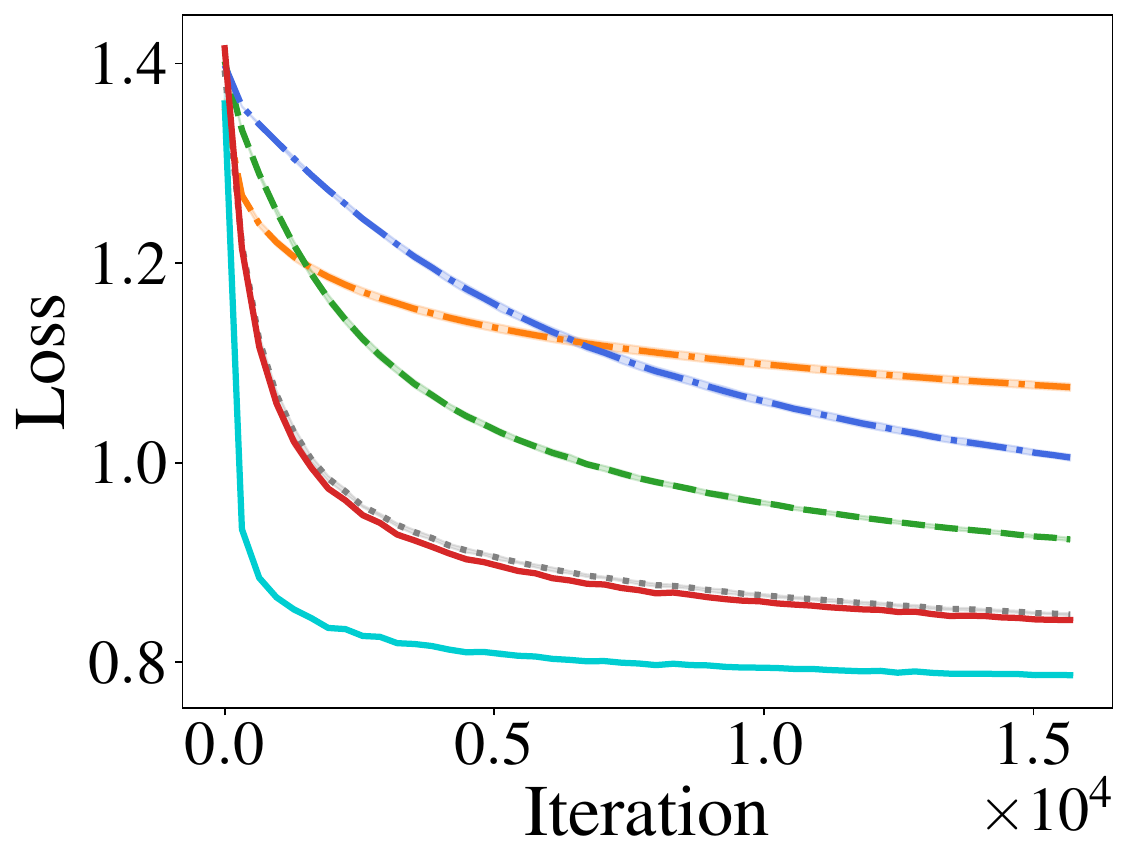}\\
\subfigure[Australian]
{\includegraphics[width=\rat\textwidth]{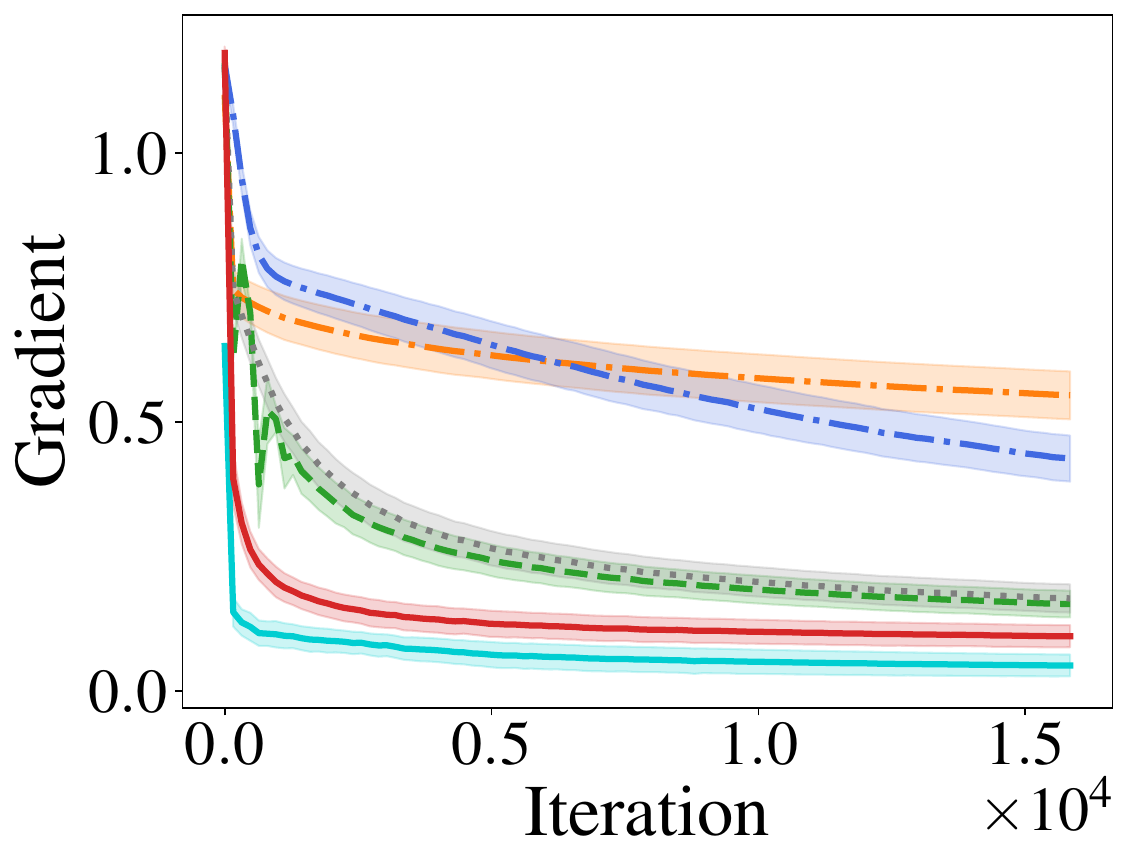}}
\subfigure[Breast-cancer]
{\includegraphics[width=\rat\textwidth]{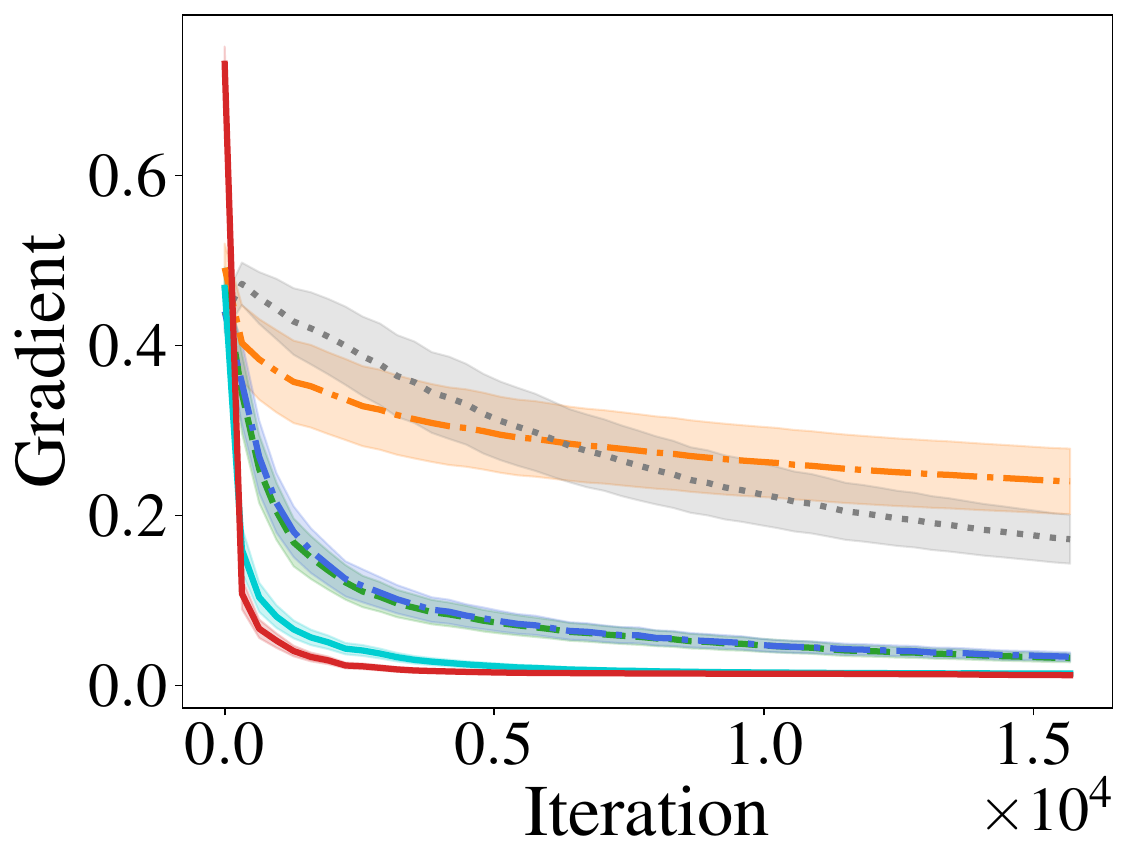}}
\subfigure[Svmguide1]
{\includegraphics[width=\rat\textwidth]{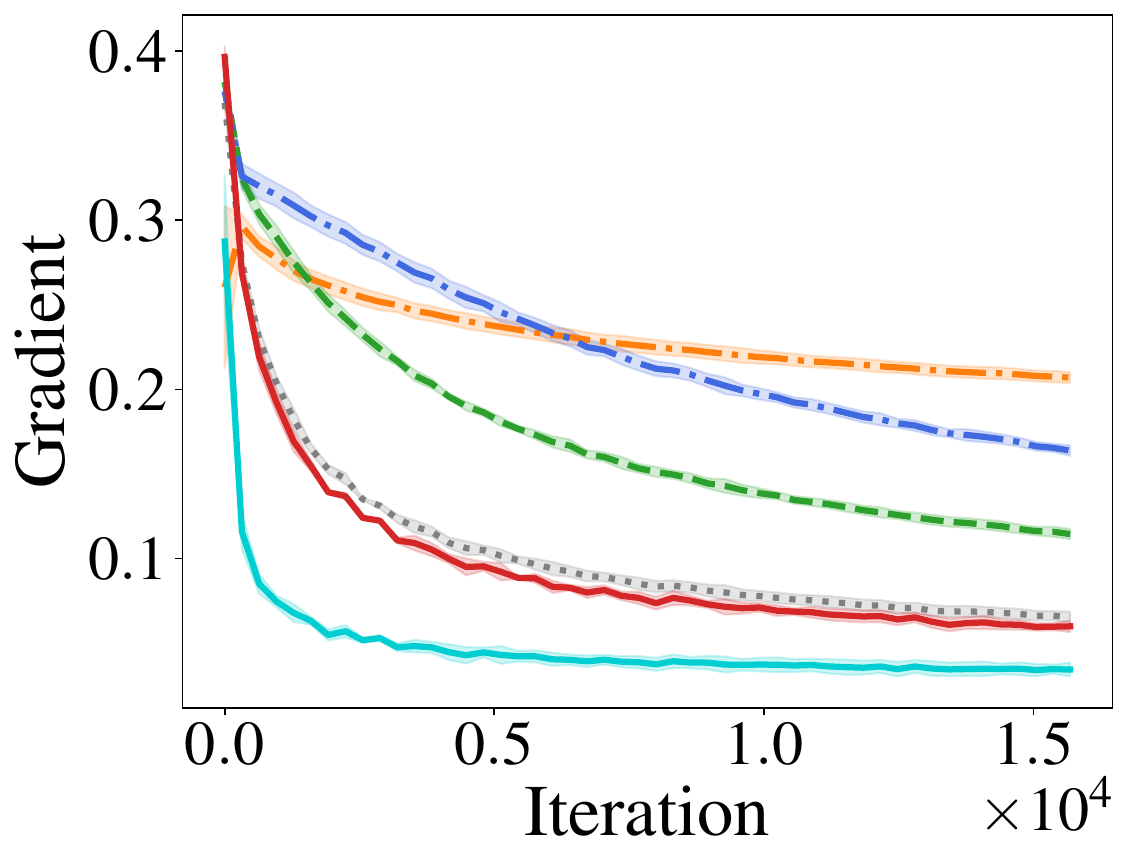}}	
{\includegraphics[width=0.9\textwidth]{./legend1.pdf}}
\caption{Results for Hierarchical Tilted Empirical Risk Minimization.}
		\label{fig:3}
	\end{center}
\end{figure*}

The Hierarchical Tilted Empirical Risk Minimization (TERM) objective~\cite{li2020tilted, li2021tilted} is a loss function designed to simultaneously achieve robustness against noisy data and ensure fairness across subgroups. 
Given samples $\{z_1,\ldots,z_N\}$ and per-sample loss $f(\theta;z_i)$, the TERM objective is defined as $\widetilde{R}(\theta):=\frac{1}{t} \log \left(\frac{1}{N} \sum_{i \in[N]} e^{t f\left(\theta;z_{i} \right)}\right)$. 
It mitigates outliers when $t<0$ and can help address class imbalance when $t>0$.
When a task involves both outliers and class imbalance, the Hierarchical TERM approach can be employed:
\begin{align*}
\tilde{J}(w):=\frac{1}{t} \log \left(\frac{1}{|\mathcal{D}|} \sum_{\mathcal{G} \subseteq \mathcal{D}}|\mathcal{G}| e^{t \tilde{R}_{\mathcal{G}}(w)}\right),
\text{ with}  \ \  \widetilde{R}_{\mathcal{G}}(w):=\frac{1}{\tau} \log \left(\frac{1}{|\mathcal{G}|} \sum_{z \in \mathcal{G}} e^{\tau \ell(w ; z)}\right),
\end{align*}
where $\mathcal{D}$ represents the set of all training samples and $\mathcal{G}$ denotes the samples for a specific class.
The parameters $t$ and $\tau$ are constants addressing outliers and class imbalance, respectively. This formulation yields a four-level stochastic compositional optimization problem:
\begin{align*}
     f_{1}(w)=\frac{1}{|\mathcal{G}|} \sum_{z \in \mathcal{G}} e^{\tau \ell(w ; z)},\  f_{2}\left(x\right)=\frac{1}{\tau}\log(x), 
     f_{3}\left(y\right)=\frac{1}{|\mathcal{D}|} \sum_{\mathcal{G}  \subseteq\mathcal{D}}|\mathcal{G}| e^{t y}, \ f_{4}\left(z\right)=\frac{1}{t}\log(z).
\end{align*}
In the experiments, we evaluate the different methods using the \texttt{australian}, \texttt{breast-cancer}, and \texttt{svmguide1} datasets from LIBSVM~\cite{CC01a}.
We inject both label noise and class imbalance: approximately 30\% of labels are randomly reshuffled, and the rare-to-common class ratio is set to 1:20.
As shown in Figure~\ref{fig:3}, MML and MMLFA decrease both the objective and gradient norm faster than the baselines.
Table~\ref{sample-table} reports classification accuracy (mean$\pm$std over 20 runs), showing that our methods achieve the best performance on both the rare class and overall accuracy.
Finally, Figure~\ref{fig:4} in Appendix~\ref{addExp} provides a sensitivity analysis for the clipping threshold $C$, indicating that performance is robust when $C$ is sufficiently large (e.g., $C>1$).
\begin{table*}[!t]
	\caption{Classification accuracy (mean$\pm$std) for Hierarchical Tilted Empirical Risk Minimization.}
	\label{sample-table}
	\centering
    \resizebox{\textwidth}{!}{
    \begin{tabular}{l|cc|cc|cc}
    \toprule
    \multirow{2}*{Method} & \multicolumn{2}{c|}{Australian} &\multicolumn{2}{c|}{Breast-cancer}  &\multicolumn{2}{c}{Svmguide1}\\
    \cline{2-7}				 &\text{rare class}&\text{overall}&\text{rare class}&\text{overall}&\text{rare class}&\text{overall}\\
    \hline
					A-TSCGD    	& 0.609$\pm$0.011&  0.717$\pm$0.007 & 0.610$\pm$0.007 & 0.819$\pm$0.001 & 0.744$\pm$0.002 & 0.861$\pm$0.001\\
					SCSC 		& 0.723$\pm$0.004&  0.782$\pm$0.002 & 0.668$\pm$0.009 & 0.841$\pm$0.004 & 0.777$\pm$0.008 & 0.881$\pm$0.003\\
					Nested-SPIDER    	& 0.692$\pm$0.010&  0.754$\pm$0.008 & 0.747$\pm$0.010 & 0.877$\pm$0.005 & 0.762$\pm$0.003 & 0.864$\pm$0.005\\
					NLASG    	& 0.750$\pm$0.006&  0.790$\pm$0.005 & 0.728$\pm$0.011 & 0.870$\pm$0.006 & 0.754$\pm$0.002 & 0.866$\pm$0.002\\
					\hline     
     \textbf{MML}    & 0.813$\pm$0.004& \textbf{0.823$\pm$0.006}& 0.805$\pm$0.004 & 0.913$\pm$0.003 & \textbf{0.782$\pm$0.002} & \textbf{0.884$\pm$0.002}\\
     \textbf{MMLFA}    & \textbf{0.838$\pm$0.007}& \text{0.817$\pm$0.004}& \textbf{0.824$\pm$0.002}& \textbf{0.920$\pm$0.002} & \text{0.778$\pm$0.001}&  \text{0.882$\pm$0.001}\\
					\bottomrule
				  \end{tabular}
            }  
\end{table*}
\section{Conclusion}
In this paper, we addressed stochastic multi-level optimization by employing momentum-based estimators to track function values and gradients.
With mini-batches, the resulting method achieves the optimal sample complexity $\mathcal{O}(\epsilon^{-4})$ under standard smoothness assumptions, without requiring the stronger average smoothness conditions used in prior work.
We further introduced a normalized variant that preserves the same rate while avoiding dependence on problem-specific constants to set hyperparameters.
To eliminate large batches, we developed a batch-free approach based on first-order approximation and clipping, which retains the optimal sample complexity with a constant batch size.
Experiments on risk-averse portfolio optimization and hierarchical tilted empirical risk minimization further demonstrate the effectiveness of the proposed methods.

\bibliography{example_paper}
\bibliographystyle{iclr2027_conference}
\newpage
\appendix
\section{Proof of Lemma~\ref{L1}}
Considering that $\u^k_i = (1-\alpha)\u^{k-1}_i  + \alpha \frac{1}{B_1}\sum_{j=1}^{B_1}f_i(\u^k_{i-1};\xi^k_{i,j}) $, we have:
\begin{align*}
& \E\left[\Norm{\u^k_i - f_i(\u^k_{i-1})}^2\right]\\
=& \E\left[\left\|(1-\alpha)\left(\u^{k-1}_i - f_i(\u^{k}_{i-1})\right)   +\alpha\left(\frac{1}{B_1}\sum_{j=1}^{B_1}   f_i(\u^{k}_{i-1};\xi^k_{i,j}) - f_i(\u^k_{i-1}) \right)\right\|^2\right] \\
=& (1-\alpha)^2 \E\left[\left\|\u^{k-1}_i - f_i(\u^{k}_{i-1})  \right\|^2\right] + \alpha^2 \E\left[\left\|\frac{1}{B_1}\sum_{j=1}^{B_1}   f_i(\u^{k}_{i-1};\xi^k_{i,j}) - f_i(\u^k_{i-1}) \right\|^2\right]\\
\leq & (1-\alpha)^2(1+\alpha) \E\left[\left\|\u^{k-1}_i - f_i(\u^{k-1}_{i-1})  \right\|^2\right]+(1-\alpha)^2(1+\frac{1}{\alpha}) \E\left[\left\|f_i(\u^{k-1}_{i-1}) - f_i(\u^{k}_{i-1})  \right\|^2\right] + \frac{\alpha^2\sigma^2}{B_1} \\
\leq & (1-\alpha) \E\left[\left\|\u^{k-1}_i - f_i(\u^{k-1}_{i-1})  \right\|^2\right]+\frac{2}{\alpha} \E\left[\left\|f_i(\u^{k-1}_{i-1}) - f_i(\u^{k}_{i-1})  \right\|^2\right] + \frac{\alpha^2\sigma^2}{B_1} \\
\leq & (1-\alpha) \E\left[\left\|\u^{k-1}_i - f_i(\u^{k-1}_{i-1})  \right\|^2\right]+\frac{2\LL_f^2}{\alpha} \E\left[\left\|\u^{k-1}_{i-1} - \u^{k}_{i-1}  \right\|^2\right] + \frac{\alpha^2\sigma^2}{B_1}.
\end{align*}
Set $B_1 \leq \alpha^2 B_0 K$. Summing up and rearranging, we can obtain:
\begin{align*}
    \frac{1}{K}\sum_{k=1}^K  \E\left[\Norm{\u^k_i -  f_i(\u^k_{i-1})}^2\right] \leq &\frac{\E\left[\Norm{\u^1_i -  f_i(\u^1_{i-1})}^2\right]}{\alpha K} + \frac{2\LL_f^2}{\alpha^2 K}  \sum_{k=1}^K \E\left[\Norm{\u^{k+1}_{i-1} - \u^{k}_{i-1}}^2\right] + \frac{\alpha   \sigma^2}{B_1} \\
 \leq &\frac{\sigma^2}{\alpha K B_0} + \frac{2\LL_f^2}{\alpha^2 K}  \sum_{k=1}^K \E\left[\Norm{\u^{k+1}_{i-1} - \u^{k}_{i-1}}^2\right] + \frac{\alpha   \sigma^2}{B_1}\\
 \leq &\frac{2\LL_f^2}{\alpha^2 K}  \sum_{k=1}^K \E\left[\Norm{\u^{k+1}_{i-1} - \u^{k}_{i-1}}^2\right] + \frac{2\alpha   \sigma^2}{B_1}.
\end{align*}
Next, we bound the difference of function estimators between adjacent steps:
\begin{align*}
        & \E\left[\Norm{\u_{i}^{k+1} - \u_{i}^{k}}^2\right] =\alpha^2\E\left[\left\|\frac{1}{B_1}\sum_{j=1}^{B_1}f_{i}(\u^{k+1}_{i-1};\xi^{k+1}_{i,j}) - \u_{i}^{k} \right\|^2 \right]\\
  =&\alpha^2\E\left[\left\|\frac{1}{B_1}\sum_{j=1}^{B_1}f_{i}(\u^{k+1}_{i-1};\xi^{k+1}_{i,j})-f_{i}(\u^{k+1}_{i-1}) + f_{i}(\u^{k+1}_{i-1}) - \u_{i}^{k} \right\|^2 \right]\\
  =&\alpha^2\E\left[\left\|\frac{1}{B_1}\sum_{j=1}^{B_1}f_{i}(\u^{k+1}_{i-1};\xi^{k+1}_{i,j})-f_{i}(\u^{k+1}_{i-1})  \right\|^2 \right] + \alpha^2\E\left[\left\| f_{i}(\u^{k+1}_{i-1}) - \u_{i}^{k} \right\|^2 \right]\\
  \leq&\frac{\alpha^2\sigma^2}{B_1}  +\alpha^2\E\left[\left\| f_{i}(\u^{k+1}_{i-1}) -f_{i}(\u^{k}_{i-1}) + f_{i}(\u^{k}_{i-1}) - \u_{i}^{k} \right\|^2 \right]\\
  \leq &\frac{\alpha^2\sigma^2}{B_1}  +2\alpha^2\E\left[\left\| f_{i}(\u^{k+1}_{i-1}) -f_{i}(\u^{k}_{i-1}) \right\|^2 \right] + 2\alpha^2\E\left[\left\| f_{i}(\u^{k}_{i-1}) - \u_{i}^{k} \right\|^2 \right]\\
  \leq &\frac{\alpha^2\sigma^2}{B_1}  +2\alpha^2\LL_f^2 \E\left[\left\| \u^{k+1}_{i-1} -\u^{k}_{i-1} \right\|^2 \right] + 2\alpha^2\E\left[\left\| f_{i}(\u^{k}_{i-1}) - \u_{i}^{k} \right\|^2 \right].
    \end{align*} 
Summing up, we obtain that
\begin{align*}
        \frac{1}{K}\sum_{k=1}^K \E\left[\Norm{\u_{i}^{k+1} - \u_{i}^{k}}^2\right] \leq & \frac{\alpha^2\sigma^2}{B_1}  + \frac{2\alpha^2\LL_f^2}{K}\sum_{k=1}^K \E\left[\left\| \u^{k+1}_{i-1} -\u^{k}_{i-1} \right\|^2 \right] +  \frac{2\alpha^2}{K}\sum_{k=1}^K \E\left[\left\| f_{i}(\u^{k}_{i-1}) - \u_{i}^{k} \right\|^2 \right].
    \end{align*} 

\section{Proof of Lemma~\ref{L2}}
We would prove this Lemma by Induction.

(1) For level $i=1$, by noting that $\u^k_0 = \x^k$ and setting $C_1=5+6\LL_f^2$, we have
\begin{align*}
    &\frac{1}{K}\sum_{k=1}^K  \E\left[\Norm{\u_1^k -  f_1(\u^k_{0})}^2\right]\\ 
 \leq &\frac{2\LL_f^2}{\alpha^2 K}  \sum_{k=1}^K \E\left[\Norm{\u^{k+1}_{0} - \u^{k}_{0}}^2\right] + \frac{2\alpha   \sigma^2}{B_1}
 \leq  \frac{2\LL_f^2}{\alpha^2 K}  \sum_{k=1}^K \E\left[\Norm{\x^{k+1} - \x^{k}}^2\right] + \frac{2\alpha   \sigma^2}{B_1}\\
 \leq & (2\LL_f^2+2) \left(\frac{1}{\alpha^2 K}  \sum_{k=1}^K \E\left[\Norm{\x^{k+1} - \x^{k}}^2\right]  + \frac{ \sigma^2}{B_1} \right) \\
 \leq & C_1 \left(\frac{1}{\alpha^2 K}  \sum_{k=1}^K \E\left[\Norm{\x^{k+1} - \x^{k}}^2\right]  + \frac{ \sigma^2}{B_1} \right),
\end{align*}
as well as the following guarantee
\begin{align*}
 &\frac{1}{K}\sum_{k=1}^K \E\left[\Norm{\u^{k+1}_{1} - \u^{k}_{1}}^2\right]\\ \leq & \frac{\alpha^2\sigma^2}{B_1}  + \frac{2\alpha^2\LL_f^2}{K}\sum_{k=1}^K \E\left[\left\| \u^{k+1}_{0} -\u^{k}_{0} \right\|^2 \right] +  \frac{2\alpha^2}{K}\sum_{k=1}^K \E\left[\left\| f_{1}(\u^{k}_{0}) - \u_{1}^{k} \right\|^2 \right]\\
 \leq & \frac{\alpha^2\sigma^2}{B_1}  + 2\alpha^2\LL_f^2  \frac{1}{K}  \sum_{k=1}^K \E\left[\Norm{\x^{k+1} - \x^{k}}^2\right]+  2\alpha^2\left(\frac{2\LL_f^2}{\alpha^2} \frac{1}{K}  \sum_{k=1}^K \E\left[\Norm{\x^{k+1} - \x^{k}}^2\right] + \frac{2\alpha   \sigma^2}{B_1}\right)\\
 \leq & (6\LL_f^2+5)\alpha^2\left(\frac{1}{\alpha^2 K}  \sum_{k=1}^K \E\left[\Norm{\x^{k+1} - \x^{k}}^2\right]  + \frac{\sigma^2}{B_1}\right) \\
 \leq & C_1 \alpha^2 \left( \frac{1}{\alpha^2 K}  \sum_{k=1}^K \E\left[\Norm{\x^{k+1} - \x^{k}}^2\right] + \frac{ \sigma^2}{B_1} \right).
\end{align*}
(2) Suppose that for level $i=n$, we have that
\begin{align*}
\frac{1}{K}\sum_{k=1}^K  \E\left[\Norm{\u^k_n -  f_n(\u^k_{n-1})}^2\right] \leq C_n \left(\frac{1}{\alpha^2 K}  \sum_{k=1}^K \E\left[\Norm{\x^{k+1} - \x^{k}}^2\right]  + \frac{ \sigma^2}{B_1} \right),\\ 
\frac{1}{K}\sum_{k=1}^K \E\left[\Norm{\u_{n}^{k+1} - \u_{n}^{k}}^2\right]\leq C_n \alpha^2 \left(\frac{1}{\alpha^2 K}  \sum_{k=1}^K \E\left[\Norm{\x^{k+1} - \x^{k}}^2\right]  + \frac{ \sigma^2}{B_1} \right).
\end{align*}
Then, for level $i=n+1$, we can prove that:
\begin{align*}
  \frac{1}{K}\sum_{k=1}^K  \E\left[\Norm{\u^k_{n+1} -  f_{n+1}(\u_n^{k})}^2\right]
  \leq & \frac{2\LL_f^2}{\alpha^2 K}  \sum_{k=1}^K \E\left[\Norm{\u_{n}^{k+1} - \u_{n}^{k}}^2\right] + \frac{2\alpha   \sigma^2}{B_1} \\
  \leq &\frac{2\LL_f^2}{\alpha^2}C_n \alpha^2 \left(\frac{1}{\alpha^2 K}  \sum_{k=1}^K \E\left[\Norm{\x^{k+1} - \x^{k}}^2\right]  + \frac{ \sigma^2}{B_1} \right) + \frac{2\alpha   \sigma^2}{B_1}\\
  \leq &(2+2\LL_f^2 C_n) \left(\frac{1}{\alpha^2 K}  \sum_{k=1}^K \E\left[\Norm{\x^{k+1} - \x^{k}}^2\right]  + \frac{ \sigma^2}{B_1} \right) \\
  \leq & C_{n+1} \left(\frac{1}{\alpha^2 K}  \sum_{k=1}^K \E\left[\Norm{\x^{k+1} - \x^{k}}^2\right]  + \frac{ \sigma^2}{B_1} \right), 
\end{align*}
as well as the fact that
\begin{align*}
 \frac{1}{K}\sum_{k=1}^K \E\left[\Norm{\u_{n+1}^{k+1} - \u_{n+1}^{k}}^2\right] 
 \leq & \frac{\alpha^2\sigma^2}{B_1}  + \frac{2\alpha^2\LL_f^2}{K}\sum_{k=1}^K \E\left[\left\| \u_{n}^{k+1} -\u_{n}^{k} \right\|^2 \right] +  \frac{2\alpha^2}{K}\sum_{k=1}^K \E\left[\left\| \u_{n+1}^{k} - f_{n+1}(\u_{n}^{k}) \right\|^2 \right]\\ 
 \leq & \frac{\alpha^2\sigma^2}{B_1}  + 2\alpha^2\LL_f^2 C_n \alpha^2 \left(\frac{1}{\alpha^2 K}  \sum_{k=1}^K \E\left[\Norm{\x^{k+1} - \x^{k}}^2\right]  + \frac{ \sigma^2}{B_1} \right) \\
 & +  2\alpha^2(2+2\LL_f^2 C_n) \left(\frac{1}{\alpha^2 K}  \sum_{k=1}^K \E\left[\Norm{\x^{k+1} - \x^{k}}^2\right]  + \frac{ \sigma^2}{B_1} \right)\\ 
 \leq & (5+6\LL_f^2 C_n) \alpha^2\left(\frac{1}{\alpha^2 K}  \sum_{k=1}^K \E\left[\Norm{\x^{k+1} - \x^{k}}^2\right]  + \frac{ \sigma^2}{B_1} \right)\\ 
 \leq & C_{n+1} \alpha^2\left(\frac{1}{\alpha^2 K}  \sum_{k=1}^K \E\left[\Norm{\x^{k+1} - \x^{k}}^2\right]  + \frac{ \sigma^2}{B_1} \right).
\end{align*}
Thus, we have proved this lemma via Induction. Also, by setting that $C_0=1$, we have $C_{i+1} = 5+6L_f^2 C_i$ for $i\ge 1$.

\section{Proof of Lemma~\ref{L3}}
First, we define $\y^k_i = f_i\circ f_{i-1}\circ \dotsc \circ f_1(\x^k)$, and define its gradient as $\nabla F_i (\x^k)$. Also denote that $\nabla \widehat{F}_i(\x^k) \coloneqq \nabla f_1(\x^k) \cdots \nabla f_i(\u_{i-1}^k)$. Then, we have the following guarantee:
\begin{align*}
        \Norm{\y^k_1 - \u_1^k} &=\Norm{ f_1(\x^k) - \u^k_1} \\
		\Norm{\y^k_2 - \u^k_2} &=\Norm{f_2\circ f_1(\x^k) - \u_2^k} \\
        &\leq \Norm{f_2\circ f_1(\x^k) - f_2(\u^k_1)} + \Norm{f_2(\u^k_1)- \u_2^k}\\
		& \leq L_f \Norm{f_1(\x^k) - \u^k_1} + \Norm{f_2(\u^k_1)- \u_2^k},\\
		\Norm{\y^k_3-\u_3^k} &= \Norm{f_3\circ f_2\circ f_1(\x^k) - \u_3^k} \\
        &\leq \Norm{f_3\circ f_2 \circ f_1(\x^k) - f_3(\u_2^k)} + \Norm{f_3(\u_2^k) - \u_3^k}\\
        & \leq L_f \Norm{f_2 \circ f_1(\x^k) - \u_2^k}  + \Norm{f_3(\u_2^k) - \u_3^k}\\
		& \leq L_f \Norm{\y^k_2 - \u_2^k}  + \Norm{f_3(\u_2^k) - \u_3^k}\\
		& \leq L_f(L_f \Norm{f_1(\x^k) - \u^k_1} + \Norm{f_2(\u_1^k) - \u^k_2}) +  \Norm{f_3(\u_2^k) - \u_3^k}\\
		& \cdots\\
		\Norm{\y^k_i-\u^k_i} & \leq L_f \Norm{\y^k_{i-1} - \u^k_{i-1}} + \Norm{f_i(\u^k_{i-1}) - \u^k_i}\\
        & \leq L_f \left(L_f \Norm{\y^k_{i-2} - \u^k_{i-2}} + \Norm{f_{i-1}(\u^k_{i-2}) - \u^k_{i-1}} \right)+ \Norm{f_i(\u^k_{i-1}) - \u^k_i}\\
		& \leq \sum_{j=1}^i L_f^{i-j} \Norm{f_j(\u^k_{j-1}) - \u^k_j}.
	\end{align*}
 Besides, we also know that:
\begin{align*}
		\Norm{\nabla F_1(\x^k) - \nabla \widehat{F}_1(\x^k)} &= 0,\\
		\Norm{\nabla F_2(\x^k) - \nabla \widehat{F}_2(\x^k)} &= \Norm{\nabla f_1(\x^k) \nabla f_2(\y_1^k) - \nabla f_1(\x^k) \nabla f_2(\u_1^k)}\leq  L_f L_J\Norm{\y^k_1 - \u_1^k},	\\
		\Norm{\nabla F_3(\x^k) - \nabla \widehat{F}_3(\x^k)} &= \Norm{\nabla f_1(\x^k) \nabla f_2(\y_1^k) \nabla f_3(\y_2^k) -  \nabla f_1(\x^k) \nabla f_2(\u^k_1) \nabla f_3(\u_2^k)}\\
        & \leq L_f \Norm{\nabla f_2(\y_1^k) \nabla f_3(\y_2^k) -   \nabla f_2(\u^k_1) \nabla f_3(\u_2^k)}\\
        & = L_f \left\|\nabla f_2(\y_1^k) \left(\nabla f_3(\y_2^k) - \nabla f_3(\u_2^k)\right)+\left( \nabla f_2(\y_1^k) - \nabla f_2(\u^k_1) \right) \nabla f_3(\u_2^k)\right\|\\
		&\leq L_f^2 L_J\left(\Norm{\y^k_2 - \u^k_2}  +  \Norm{\y^k_1 - \u^k_1}\right),\\
		& \cdots
	\end{align*}
	To this end, we can conclude that:
	\begin{align*}
		\Norm{\nabla F_T(\x^k) - \nabla\widehat{F}_T(\x^k)}  &\leq L_f^{T-1}L_J \sum_{i=1}^{T-1} \Norm{\y^k_i - \u^k_i}\\
		 & \leq L_f^{T-1} L_J(1+L_f+\dotsc+L_f^{T-1}) \sum_{i=1}^{T-1}  \Norm{f_i(\u_{i-1}^k) - \u_i^k}\leq L_F \sum_{i=1}^{T-1}  \Norm{f_i(\u_{i-1}^k) - \u_i^k},
	\end{align*}
  and thus we can obtain
  \begin{align*}
		\Norm{\nabla F(\x^k) -\prod_{i=1}^T\nabla f_i(\u^k_{i-1})}^2 & \leq L_F^2 T  \sum_{i=1}^{T-1}  \Norm{f_i(\u_{i-1}^k) - \u_i^k}^2.
	\end{align*}
Next, considering the definition of $\z^{k+1}$, we have
    \begin{align*}
        &\E\left[\left\|\z^{k+1} - \nabla F\left(\x^{k+1}\right) \right\|^2\right] \\
        = & \E\left[\Norm{(1-\alpha)(\z^{k}-\nabla F(\x^{k+1}))  +\alpha\left( \prod_{i=1}^T \nabla f_i(\u^{k+1}_{i-1}) -\nabla F(\x^{k+1})\right)\right.\right.\\
        &\quad \left.\left.+\alpha\left( \frac{1}{B_1} \sum_{j=1}^{B_1} \left[ \prod_{i=1}^T \nabla f_i(\u^{k+1}_{i-1};\xi^{k+1}_{i,j}) \right] -\prod_{i=1}^T \nabla f_i(\u^{k+1}_{i-1})\right)}^2\right]\\
        \leq & \E\left[\Norm{(1-\alpha)(\z^{k}-\nabla F(\x^{k+1}))  +\alpha\left( \prod_{i=1}^T \nabla f_i(\u^{k+1}_{i-1}) -\nabla F(\x^{k+1})\right)}^2 \right]\\
        &+\alpha^2\E\left[\Norm{\frac{1}{B_1} \sum_{j=1}^{B_1} \left[ \prod_{i=1}^T \nabla f_i(\u^{k+1}_{i-1};\xi^{k+1}_{i,j}) \right] -\prod_{i=1}^T \nabla f_i(\u^{k+1}_{i-1})}^2\right]\\
        \leq & \E\left[\Norm{(1-\alpha)(\z^{k}-\nabla F(\x^{k}))+(1-\alpha)(\nabla F(\x^{k}) -\nabla F(\x^{k+1}))  +\alpha\left( \prod_{i=1}^T \nabla f_i(\u^{k+1}_{i-1}) -\nabla F(\x^{k+1})\right)}^2\right] \\
        &+\frac{\alpha^2}{B_1}T^2 (\sigma^2 +L_f^2 )^{T-1}\sigma^2\\
        \leq & (1-\alpha)\E\left[\Norm{\z^{k}-\nabla F(\x^{k})}^2\right] +\frac{4}{\alpha}\E\left[\Norm{\nabla F(\x^{k}) -\nabla F(\x^{k+1})}^2\right]+4\alpha\E\left[\Norm{ \prod_{i=1}^T \nabla f_i(\u^{k+1}_{i-1}) -\nabla F(\x^{k+1})}^2\right]\\
        &+\frac{\alpha^2}{B_1}T^2 (\sigma^2 +L_f^2 )^{T-1}\sigma^2\\
        \leq & (1-\alpha)\E\left[\Norm{\z^{k}-\nabla F(\x^{k})}^2\right] +\frac{4L_F^2}{\alpha}   \E\left[\Norm{\x^{k+1} - \x^{k}}^2\right] +4\alpha T L_F^2 \sum_{i=1}^{T-1} \E\left[\Norm{f_i(\u^{k+1}_{i-1})-\u^{k+1}_i}^2\right]\\
        &+\frac{\alpha^2}{B_1}T^2 (\sigma^2 +L_f^2 )^{T}.
    \end{align*}
Summing up and rearranging, we proved the conclusion by noting that $B_1 \leq \alpha^2 B_0 K$ and
 \begin{align*}
     \E\left[\Norm{\z^{1}-\nabla F(\x^{1})}^2\right] \leq &T^2{(\sigma^2+L_f^2)^T}/{B_0}\leq \alpha^2 K T^2{(\sigma^2+L_f^2)^T}/{B_1},\\ 
     \frac{1}{K}\sum_{k=1}^K  \E\left[\Norm{\u^{k+1}_i -  f_i(\u^{k+1}_{i-1})}^2\right] &\leq \frac{1}{K}\sum_{k=1}^{K+1}  \E\left[\Norm{\u^k_i -  f_i(\u^k_{i-1})}^2\right] \leq  C_i \left(\frac{1}{\alpha^2 K}  \sum_{k=1}^{K+1} \E\left[\Norm{\x^{k+1} - \x^{k}}^2\right]  + \frac{ \sigma^2}{B_1} \right).
 \end{align*}

\section{Proof of Theorem~\ref{T1}}
Due to the smoothness of function $F$ and the definition of $\x^{k+1}$, by setting $\eta \leq \frac{1}{2L_F}$, we have:
\begin{align*}
F\left(\x^{k+1}\right) & \leq F\left(\x^{k}\right)+\left\langle\nabla F(\x^{k}), \x^{k+1}-\x^{k}\right\rangle+\frac{L_{F}}{2}\left\|\x^{k+1}-\x^{k}\right\|^{2} \\ &=F(\x^{k})-\eta\left\langle\nabla F(\x^{k}), \z^{k}\right\rangle+\frac{\eta^{2} L_{F}}{2}\|\z^{k}\|^{2} \\ & \leq F(\x^{k})-\eta\langle\nabla F(\x^{k}), \z^{k}\rangle+\frac{\eta}{2}\|\nabla F(\x^{k})\|^{2}+\frac{\eta}{2}\left\|\z^{k}\right\|^{2} -\frac{\eta}{2}\|\nabla F(\x^{k})\|^{2}-\frac{\eta}{2}\left\|\z^{k}\right\|^{2} +\frac{\eta}{4}\left\|\z^{k}\right\|^{2}\\ &=F(\x^{k})+\frac{\eta}{2}\|\nabla F(\x^{k})-\z^{k}\|^{2}-\frac{\eta}{2}\|\nabla F(\x^{k})\|^{2}-\frac{\eta}{4}\left\|\z^{k}\right\|^{2}.
\end{align*}
    Also note that $\E\left[ \Norm{\z^1}^2\right] \leq 2(L_f^2 +\sigma^2)^{T}$ and
\begin{align*}
   \E\left[ \Norm{\z^k}^2\right]
   = & \E\left[\Norm{(1-\alpha)\z^{k-1}  +\alpha \frac{1}{B_1} \sum_{j=1}^{B_1}  \prod_{i=1}^T \nabla f_i(\u^{k}_{i-1};\xi^{k}_{i,j})}^2\right]\\
\leq & (1-\alpha)^2(1+\alpha)\E\left[\Norm{\z^{k-1}  }^2\right] + \alpha^2(1+\frac{1}{\alpha})\E\left[\Norm{ \frac{1}{B_1}\sum_{j=1}^{B_1}   \prod_{i=1}^T \nabla f_i(\u^{k}_{i-1};\xi^{k}_{i,j})}^2\right]\\
\leq & (1-\alpha)\E\left[\Norm{\z^{k-1}  }^2\right]+ 2\alpha \frac{1}{B_1}\sum_{j=1}^{B_1}  \E\left[\Norm{   \prod_{i=1}^T \nabla f_i(\u^{k}_{i-1};\xi^{k}_{i,j}) }^2\right]\\
\leq & (1-\alpha)\E\left[\Norm{\z^{k-1}}^2\right] + 2\alpha \frac{1}{B_1}\sum_{j=1}^{B_1}  \prod_{i=1}^T \E_{\xi_{i,j}}\left[\Norm{    \nabla f_i(\u^{k}_{i-1};\xi^{k}_{i,j}) }^2\right]\\
\leq & (1-\alpha)\E\left[\Norm{\z^{k-1}  }^2\right] + \alpha \cdot 2 \left(L_f^2 +\sigma^2\right)^T \leq 2 \left(L_f^2 +\sigma^2\right)^T .
\end{align*}
Set $B_0 = 1$, $\alpha =\mathcal{O}(1)$, $B_1 = \alpha^2 K =\mathcal{O}(K)$ and learning rate $\eta \leq  \frac{\alpha}{4L_F \sqrt{1+T \sum_{i=1}^{T-1}C_i}} = \mathcal{O}(1)$. We can ensure
\begin{align*}
        &\E\left[\frac{1}{K}\sum_{k=1}^K \Norm{\nabla F(\x^k)}^2 \right] \\
  \leq & \frac{2\E\left[F(\x^1)-F(\x^{K+1})\right]}{\eta K} + \frac{1}{K}\sum_{k=1}^K \E\left[\Norm{\nabla F(\x^k) - \z^k}^2\right]  - \frac{1}{2K}\E\left[\sum_{k=1}^K\Norm{\z^k}^2\right] \\
  \leq & \frac{2\E\left[F(\x^1)-F(\x^{K+1})\right]}{\eta K} - \frac{1}{2K}\sum_{k=1}^K \E\left[\Norm{\z^k}^2\right]+\frac{2\alpha T^2}{B_1} (\sigma^2 +L_f^2 )^{T} +\frac{4L_F^2 \eta^2}{\alpha^2 K}  \sum_{k=1}^K \E\left[\Norm{\z^{k}}^2\right] \\
    &+4T L_F^2\sum_{i=1}^{T-1} C_i \left(\frac{\eta^2}{\alpha^2 K}  \sum_{k=1}^{K+1} \E\left[\Norm{\z^{k}}^2\right]  + \frac{ \sigma^2}{B_1} \right)\\
   \leq & \frac{2\E\left[F(\x^1)-F(\x^{K+1})\right]}{\eta K} +\frac{2\alpha T^2}{B_1} (\sigma^2 +L_f^2 )^{T} + \frac{ \sigma^2}{B_1} 4T L_F^2\sum_{i=1}^{T-1} C_i +\left(4T L_F^2\sum_{i=1}^{T-1} C_i \right)\frac{\eta^2}{\alpha^2 K}   \E\left[\Norm{\z^{K+1}}^2\right]   \\
 \leq & \mathcal{O}\left( \frac{1}{\eta K} + \frac{1}{B_1} +\frac{\eta^2}{\alpha^2 K}\right)
    \leq \mathcal{O}\left(\frac{1}{K} \right).
    \end{align*}
We complete the proof of the theorem by noticing $\E\left[\frac{1}{K}\sum_{k=1}^K \Norm{\nabla F(\x^k)} \right] \leq \sqrt{\E\left[\frac{1}{K}\sum_{k=1}^K \Norm{\nabla F(\x^k)}^2 \right]}$.

\newpage
\section{Proof of Theorem~\ref{T2}}
By noting that $\x^{k+1} = \x^k- \frac{\eta}{\Norm{\z^k}} \z^k$, we have:
\begin{align*}
F\left(\x^{k+1}\right)  \leq &F\left(\x^{k}\right)+\left\langle\nabla F(\x^{k}), \x^{k+1}-\x^{k}\right\rangle+\frac{L_{F}}{2}\left\|\x^{k+1}-\x^{k}\right\|^{2} \\ 
=&F(\x^{k})-\eta\left\langle\nabla F(\x^{k}), \frac{\z^{k}}{\Norm{\z^k}}\right\rangle+\frac{\eta^{2} L_{F}}{2} \\ 
=&F(\x^{k})   +  \eta\left\langle\z^k - \nabla F(\x^{k}), \frac{\z^{k}}{\Norm{\z^k}}\right\rangle - \eta \left\langle\z^k, \frac{\z^{k}}{\Norm{\z^k}}\right\rangle+\frac{\eta^{2} L_{F}}{2} \\ 
\leq & F(\x^k) + \eta \Norm{\z^k - \nabla F(\x^k)} \frac{\Norm{\z^k}}{\Norm{\z^k}} - \eta \frac{\Norm{\z^k}^2}{\Norm{\z^k}} + \frac{\eta^{2} L_{F}}{2}\\
\leq & F(\x^k) + \eta \Norm{\z^k - \nabla F(\x^k)} - \eta \Norm{\z^k} + \frac{\eta^{2} L_{F}}{2}.
\end{align*}   
This leads to the fact that
\begin{align*}
\Norm{\nabla F(\x^k)}\leq & \Norm{\z^k} + \Norm{\z^k - \nabla F(\x^k)}
    \leq \frac{F(\x^k) - F(\x^{k+1})}{\eta}  +  2\Norm{\z^k - \nabla F(\x^k)} +\frac{\eta L_F}{2} .
\end{align*}
Summing up, by setting $\eta = 1/\sqrt{K}$, $B_0=1$, $B_1=\lfloor K/4 \rfloor$, $\alpha = 1/2$~(actually $\alpha$ can be any positive constant within the range $(0,1)$), we know that
\begin{align*}
     &\frac{1}{K}\sum_{k=1}^K \E\left[\Norm{\nabla F(\x^k)}\right] \\
     \leq &\frac{\E\left[F(\x^{1})-F(\x^{K+1})\right]}{\eta K} + \E\left[\frac{2}{K}\sum_{k=1}^K \Norm{\z^k - \nabla F(\x^k)}\right] + \frac{\eta L_F}{2}\\
     \leq &\frac{\E\left[F(\x^{1})-F(\x^{K+1})\right]}{\eta K}  + \frac{\eta L_F}{2}\\
     &+ 2\sqrt{\frac{2\alpha T^2}{B_1} (\sigma^2 +L_f^2 )^{T} +\frac{4L_F^2}{\alpha^2 K}  \sum_{k=1}^{K+1} \E\left[\Norm{\x^{k+1} - \x^{k}}^2\right]+4T L_F^2\sum_{i=1}^{T-1} C_i \left(\frac{1}{\alpha^2 K}  \sum_{k=1}^{K+1} \E\left[\Norm{\x^{k+1} - \x^{k}}^2\right]  + \frac{ \sigma^2}{B_1} \right)}\\
     \leq & \mathcal{O}\left(\frac{1}{\eta K} + \eta +\sqrt{\frac{ \eta^2}{\alpha^2 }   +\frac{1}{B_1}} \right) \leq \mathcal{O}\left(\frac{1}{\sqrt{K}}\right).
\end{align*}

\newpage
\section{Proof of Lemma~\ref{L1+}}
Considering  $\Tilde{\u}^k_i = (1-\alpha)\left(\u^{k-1}_i+\langle \nabla f_i(\u_{i-1}^{k};\xi_{i}^k), \u_{i-1}^k-\u_{i-1}^{k-1}\rangle \right) + \alpha f_i(\u^k_{i-1};\xi^k_{i})$, $\u^k_i = \operatorname{Clip}_{C} \left[ \Tilde{\u}^k_i\right]$ and $ \Norm{f_i(\u_{i-1}^k)} \leq C$, note that this operation can be viewed as a projection operation, and with the non-expansive property, we have:
\begin{align*}
& \E\left[\Norm{\u^k_i - f_i(\u^k_{i-1})}^2\right]\\
 \leq& \E\left[\Norm{\Tilde{\u}^k_i - f_i(\u^k_{i-1})}^2\right] \\
=& \E\left[\left\|(1-\alpha)\left(\u^{k-1}_i - f_i(\u^{k}_{i-1})+\langle \nabla f_i(\u_{i-1}^{k};\xi_{i}^k), \u_{i-1}^k-\u_{i-1}^{k-1}\rangle\right)   +\alpha\left( f_i(\u^{k}_{i-1};\xi^k_{i}) - f_i(\u^k_{i-1}) \right)\right\|^2\right] \\
=& \E\left[\left\|(1-\alpha)\left(\u^{k-1}_i - f_i(\u^{k}_{i-1})+\langle \nabla f_i(\u_{i-1}^{k}), \u_{i-1}^k-\u_{i-1}^{k-1}\rangle\right)\right.\right.\\
 &\left.\left. +(1-\alpha)\langle \nabla f_i(\u_{i-1}^{k};\xi_{i}^k)-\nabla f_i(\u_{i-1}^{k}), \u_{i-1}^k-\u_{i-1}^{k-1}\rangle+\alpha\left( f_i(\u^{k}_{i-1};\xi^k_{i}) - f_i(\u^k_{i-1}) \right)\right\|^2\right] \\
=& (1-\alpha)^2 \E\left[\left\|\u^{k-1}_i - f_i(\u^{k}_{i-1})+\langle \nabla f_i(\u_{i-1}^{k}), \u_{i-1}^k-\u_{i-1}^{k-1}\rangle  \right\|^2\right] \\
&+  \E\left[\left\|(1-\alpha)\langle \nabla f_i(\u_{i-1}^{k};\xi_{i}^k)-\nabla f_i(\u_{i-1}^{k}), \u_{i-1}^k-\u_{i-1}^{k-1}\rangle+\alpha\left( f_i(\u^{k}_{i-1};\xi^k_{i}) - f_i(\u^k_{i-1}) \right) \right\|^2\right]\\
\leq & (1-\alpha)^2 \E\left[\left\|\u^{k-1}_i - f_i(\u^{k-1}_{i-1})  \right\|^2+\left\|f_i(\u^{k-1}_{i-1}) - f_i(\u^{k}_{i-1})+\langle \nabla f_i(\u_{i-1}^{k}), \u_{i-1}^k-\u_{i-1}^{k-1}\rangle  \right\|^2\right]\\
&+ (1-\alpha)^2 \E\left[2\left\|\u^{k-1}_i - f_i(\u^{k-1}_{i-1})  \right\|\left\|f_i(\u^{k-1}_{i-1}) - f_i(\u^{k}_{i-1})+\langle \nabla f_i(\u_{i-1}^{k}), \u_{i-1}^k-\u_{i-1}^{k-1}\rangle  \right\|\right]\\
&+  2(1-\alpha)^2\E\left[\left\|\langle \nabla f_i(\u_{i-1}^{k};\xi_{i}^k)-\nabla f_i(\u_{i-1}^{k}), \u_{i-1}^k-\u_{i-1}^{k-1}\rangle \right\|^2\right]+2\alpha^2\E\left[\left\| f_i(\u^{k}_{i-1};\xi^k_{i}) - f_i(\u^k_{i-1}) \right\|^2\right]\\
\leq & (1-\alpha)^2 \E\left[\left\|\u^{k-1}_i - f_i(\u^{k-1}_{i-1})  \right\|^2+2\left\|f_i(\u^{k-1}_{i-1}) - f_i(\u^{k}_{i-1})  \right\|^2+2\left\|\langle \nabla f_i(\u_{i-1}^{k}), \u_{i-1}^k-\u_{i-1}^{k-1}\rangle  \right\|^2\right]\\
& +(1-\alpha)^2 \E\left[L_J \left\|\u^{k-1}_i - f_i(\u^{k-1}_{i-1})  \right\|\left\|\u^{k-1}_{i-1} - \u^{k}_{i-1}  \right\|^2\right]\\
&+  2\E\left[\Norm{ \nabla f_i(\u_{i-1}^{k};\xi_{i}^k)-\nabla f_i(\u_{i-1}^{k})}^2 \Norm{ \u_{i-1}^k-\u_{i-1}^{k-1}}^2\right]+2\alpha^2\sigma^2\\
\leq & (1-\alpha)^2 \E\left[\left\|\u^{k-1}_i - f_i(\u^{k-1}_{i-1})  \right\|^2+2L_f^2\left\|\u^{k-1}_{i-1} - \u^{k}_{i-1}  \right\|^2+2L_f^2\left\| \u_{i-1}^k-\u_{i-1}^{k-1} \right\|^2\right]\\
& +(1-\alpha)^2 \E\left[L_J \left(\left\|\u^{k-1}_i  \right\|+\left\|f_i(\u^{k-1}_{i-1})  \right\|\right)\left\|\u^{k-1}_{i-1} - \u^{k}_{i-1}  \right\|^2\right]+  2\sigma^2 \E\left[ \Norm{ \u_{i-1}^k-\u_{i-1}^{k-1}}^2\right]+2\alpha^2\sigma^2\\
\leq & (1-\alpha)^2 \E\left[\left\|\u^{k-1}_i - f_i(\u^{k-1}_{i-1})  \right\|^2\right]+4L_f^2\E\left[\left\|\u^{k-1}_{i-1} - \u^{k}_{i-1}  \right\|^2\right]\\
& +2L_J C  \E\left[\left\|\u^{k-1}_{i-1} - \u^{k}_{i-1}  \right\|^2\right]+  2\sigma^2 \E\left[ \Norm{ \u_{i-1}^k-\u_{i-1}^{k-1}}^2\right]+2\alpha^2\sigma^2\\
\leq & (1-\alpha) \E\left[\left\|\u^{k-1}_i - f_i(\u^{k-1}_{i-1})  \right\|^2\right]+(4L_f^2+2L_JC +2\sigma^2)\E\left[\left\|\u^{k}_{i-1} - \u^{k-1}_{i-1}  \right\|^2\right]+2\alpha^2\sigma^2.
\end{align*}
Summing up and rearranging, we can obtain:
\begin{align*}
    &\frac{1}{K}\sum_{k=1}^K  \E\left[\Norm{\u^k_i -  f_i(\u^k_{i-1})}^2\right] \\
    \leq &\frac{\E\left[\Norm{\u^1_i -  f_i(\u^1_{i-1})}^2\right]}{\alpha K} + \frac{(4L_f^2+2L_JC +2\sigma^2)}{\alpha K}  \sum_{k=1}^K \E\left[\Norm{\u^{k+1}_{i-1} - \u^{k}_{i-1}}^2\right] + 2\alpha   \sigma^2 \\
 \leq &\frac{\sigma^2}{\alpha K} + \frac{(4L_f^2+2L_JC +2\sigma^2)}{\alpha K}  \sum_{k=1}^K \E\left[\Norm{\u^{k+1}_{i-1} - \u^{k}_{i-1}}^2\right] + 2\alpha   \sigma^2\\
 \leq &\frac{(4L_f^2+2L_JC +2\sigma^2)}{\alpha K}  \sum_{k=1}^K \E\left[\Norm{\u^{k+1}_{i-1} - \u^{k}_{i-1}}^2\right] + 3\alpha   \sigma^2.
\end{align*}
where the last inequality is by setting $ \alpha^2 K \ge 1$.

Next, we bound the difference of function estimators between adjacent steps:
\begin{align*}
        & \E\left[\Norm{\u_{i}^{k+1} - \u_{i}^{k}}^2\right] \leq \E\left[\Norm{\Tilde{\u}_{i}^{k+1} - \u_{i}^{k}}^2\right]\\
        \leq & 2\alpha^2\E\left[\left\|f_{i}(\u^{k+1}_{i-1};\xi^{k+1}_{i}) - \u_{i}^{k} \right\|^2 \right] + 2(1-\alpha)^2 \E\left[\Norm{\langle \nabla f_i(\u_{i-1}^{k+1};\xi_{i}^{k+1}), \u_{i-1}^{k+1}-\u_{i-1}^{k}\rangle}^2 \right]\\
  \leq &2\alpha^2\E\left[\left\|f_{i}(\u^{k+1}_{i-1};\xi^{k+1}_{i})-f_{i}(\u^{k+1}_{i-1}) + f_{i}(\u^{k+1}_{i-1}) - \u_{i}^{k} \right\|^2 \right]\\
  & +2 \E\left[\Norm{\langle \nabla f_i(\u_{i-1}^{k+1};\xi_{i}^{k+1})- \nabla f_i(\u_{i-1}^{k+1}), \u_{i-1}^{k+1}-\u_{i-1}^{k}\rangle + \langle  \nabla f_i(\u_{i-1}^{k+1}), \u_{i-1}^{k+1}-\u_{i-1}^{k}\rangle}^2 \right]\\
  =&2\alpha^2\E\left[\left\|f_{i}(\u^{k+1}_{i-1};\xi^{k+1}_{i})-f_{i}(\u^{k+1}_{i-1})  \right\|^2 \right] + 2\alpha^2\E\left[\left\| f_{i}(\u^{k+1}_{i-1}) - \u_{i}^{k} \right\|^2 \right]\\
    & +2 \E\left[\Norm{\langle \nabla f_i(\u_{i-1}^{k+1};\xi_{i}^{k+1})- \nabla f_i(\u_{i-1}^{k+1}), \u_{i-1}^{k+1}-\u_{i-1}^{k}\rangle}^2 \right] +2 \E\left[\Norm{ \langle  \nabla f_i(\u_{i-1}^{k+1}), \u_{i-1}^{k+1}-\u_{i-1}^{k}\rangle}^2 \right]\\
    \leq &{2\alpha^2\sigma^2}  +2\alpha^2\E\left[\left\| f_{i}(\u^{k+1}_{i-1}) -f_{i}(\u^{k}_{i-1}) + f_{i}(\u^{k}_{i-1}) - \u_{i}^{k} \right\|^2 \right]\\
        & +2 \E\left[\Norm{ \nabla f_i(\u_{i-1}^{k+1};\xi_{i}^{k+1})- \nabla f_i(\u_{i-1}^{k+1})}^2\Norm{ \u_{i-1}^{k+1}-\u_{i-1}^{k}}^2 \right] +2 \E\left[\Norm{  \nabla f_i(\u_{i-1}^{k+1})}^2\Norm{ \u_{i-1}^{k+1}-\u_{i-1}^{k}}^2 \right]\\
  \leq &{2\alpha^2\sigma^2}  +4\alpha^2\E\left[\left\| f_{i}(\u^{k+1}_{i-1}) -f_{i}(\u^{k}_{i-1}) \right\|^2 \right] + 4\alpha^2\E\left[\left\| f_{i}(\u^{k}_{i-1}) - \u_{i}^{k} \right\|^2 \right]\\
  & +2\sigma^2 \E\left[\Norm{ \u_{i-1}^{k+1}-\u_{i-1}^{k}}^2 \right] +2L_f^2 \E\left[\Norm{ \u_{i-1}^{k+1}-\u_{i-1}^{k}}^2 \right]\\
  \leq &{2\alpha^2\sigma^2} +(2\sigma^2+6\LL_f^2) \E\left[\left\| \u^{k+1}_{i-1} -\u^{k}_{i-1} \right\|^2 \right] + 4\alpha^2\E\left[\left\| f_{i}(\u^{k}_{i-1}) - \u_{i}^{k} \right\|^2 \right].
    \end{align*} 
Summing up, we obtain that
\begin{align*}
        \frac{1}{K}\sum_{k=1}^K \E\left[\Norm{\u_{i}^{k+1} - \u_{i}^{k}}^2\right] \leq & 2\alpha^2\sigma^2  + \frac{6\LL_f^2+2\sigma^2}{K}\sum_{k=1}^K \E\left[\left\| \u^{k+1}_{i-1} -\u^{k}_{i-1} \right\|^2 \right] +  \frac{4\alpha^2}{K}\sum_{k=1}^K \E\left[\left\| f_{i}(\u^{k}_{i-1}) - \u_{i}^{k} \right\|^2 \right].
    \end{align*}

\section{Proof of Lemma~\ref{L2+}}
Note that we have proved that
\begin{align*}
	    \Gamma_i \leq \frac{(4L_f^2+2L_JC +2\sigma^2)}{\alpha}  \Lambda_{i-1} + 3\alpha   \sigma^2; \quad 
 \Lambda_i \leq  2\alpha^2\sigma^2  + (6\LL_f^2+2\sigma^2)\Lambda_{i-1} +  (4\alpha^2)\Gamma_i.
\end{align*}
Thus by setting $\alpha \leq 1/2$, we know that
\begin{align*}
	    \Gamma_i +\frac{\Lambda_i}{2\alpha}\leq& \frac{(4L_f^2+2L_JC +2\sigma^2)}{\alpha}  \Lambda_{i-1} + 3\alpha   \sigma^2+\alpha\sigma^2  + \frac{(3\LL_f^2+\sigma^2)\Lambda_{i-1}}{\alpha} +  (2\alpha)\Gamma_i\\
     \leq &\frac{(7L_f^2+2L_JC +3\sigma^2)}{\alpha}  \Lambda_{i-1} + 4\alpha   \sigma^2 +  \Gamma_i.
\end{align*}
As a result, denoting $D = 2(7L_f^2+2L_JC +3\sigma^2)$ and noting $\Lambda_0=\frac{1}{K}\sum_{k=1}^K \E\left[\Norm{\x^{k+1} - \x^{k}}^2\right]$, we obtain
\begin{align*}
	    \Lambda_i 
     \leq &D  \Lambda_{i-1} + 8\alpha^2   \sigma^2
     \leq D^2  \Lambda_{i-2} + 8\alpha^2   \sigma^2(1+D)
     \leq D^i  \Lambda_{0} + 8\alpha^2   \sigma^2 D^i \\
     \leq & D^i \frac{1}{K}\sum_{k=1}^K\E\left[\Norm{\x^{k+1} - \x^{k}}^2\right] + 8\alpha^2   \sigma^2 D^i.
\end{align*}
Then, we can obtain that
\begin{align*}
	    \Gamma_i \leq & \frac{D}{\alpha}  \left[D^{i-1} \frac{1}{K}\sum_{k=1}^K\E\left[\Norm{\x^{k+1} - \x^{k}}^2\right] + 8\alpha^2   \sigma^2 D^{i-1}\right] + 3\alpha   \sigma^2\\
     \leq & \frac{D^i}{\alpha} \frac{1}{K}\sum_{k=1}^K\E\left[\Norm{\x^{k+1} - \x^{k}}^2\right]  + (8 D^{i} +3)\alpha   \sigma^2.
\end{align*}

\newpage
\section{Proof of Lemma~\ref{L3+}}
In the previous analysis, we have already proved that
    \begin{align*}
        \E\left[\left\|\z^{k+1} - \nabla F\left(\x^{k+1}\right) \right\|^2\right] 
        \leq (1-\alpha)\E\left[\Norm{\z^{k}-\nabla F(\x^{k})}^2\right] +\frac{4L_F^2}{\alpha}   \E\left[\Norm{\x^{k+1} - \x^{k}}^2\right] \\
        +4\alpha T L_F^2 \sum_{i=1}^{T-1} \E\left[\Norm{f_i(\u^{k+1}_{i-1})-\u^{k+1}_i}^2\right]+\alpha^2{T^2} (\sigma^2 +L_f^2 )^{T}.
    \end{align*}
Summing up and rearranging, 
    \begin{align*}
        &\frac{1}{K}\sum_{k=1}^K \E\left[\Norm{\z^{k}-\nabla F(\x^{k})}^2\right]& \\
        \leq & \frac{\E\left[\Norm{\z^{1}-\nabla F(\x^{1})}^2\right]}{\alpha K} +\frac{4L_F^2}{\alpha^2 K}\sum_{k=1}^K   \E\left[\Norm{\x^{k+1} - \x^{k}}^2\right] + \frac{4 T L_F^2}{K}\sum_{k=1}^K\sum_{i=1}^{T-1} \E\left[\Norm{f_i(\u^{k+1}_{i-1})-\u^{k+1}_i}^2\right]+\alpha {T^2} (\sigma^2 +L_f^2 )^{T} \\
\leq & \alpha T^2{(\sigma^2+L_f^2)^T}+\frac{4L_F^2}{\alpha^2 K}\sum_{k=1}^K   \E\left[\Norm{\x^{k+1} - \x^{k}}^2\right] + \frac{4 T L_F^2}{K}\sum_{k=1}^{K+1}\sum_{i=1}^{T-1} \E\left[\Norm{f_i(\u^{k}_{i-1})-\u^{k}_i}^2\right]+\alpha {T^2} (\sigma^2 +L_f^2 )^{T} \\
        \leq & 2\alpha T^2{(\sigma^2+L_f^2)^T} +\frac{4L_F^2}{\alpha^2 K}\sum_{k=1}^K   \E\left[\Norm{\x^{k+1} - \x^{k}}^2\right] + (4 T L_F^2)\sum_{i=1}^{T-1}\left( \frac{D^i}{\alpha} \frac{1}{K}\sum_{k=1}^{K+1}\E\left[\Norm{\x^{k+1} - \x^{k}}^2\right]  + (8 D^{i} +3)\alpha   \sigma^2\right)\\
 \leq & 2\alpha T^2{(\sigma^2+L_f^2)^T} +\frac{4L_F^2}{\alpha^2 K}\sum_{k=1}^K   \E\left[\Norm{\x^{k+1} - \x^{k}}^2\right] + \frac{4 T L_F^2 D^T}{\alpha} \frac{1}{K}\sum_{k=1}^{K+1}\E\left[\Norm{\x^{k+1} - \x^{k}}^2\right] +4 T L_F^2   \sigma^2(8 D^{T} +3T) \alpha \\
 \leq & \left(2 T^2{(\sigma^2+L_f^2)^T} +4 T L_F^2   \sigma^2(8 D^{T} +3T) \right)\alpha + 4L_F^2\left(1+ T  D^T\right) \frac{1}{\alpha^2K}\sum_{k=1}^{K+1}\E\left[\Norm{\x^{k+1} - \x^{k}}^2\right]\\
  \leq & D_0\alpha + 4L_F^2\left(1+ T  D^T\right) \frac{1}{\alpha^2K}\sum_{k=1}^{K+1}\E\left[\Norm{\x^{k+1} - \x^{k}}^2\right],
    \end{align*}
where $D_0 = 2 T^2{(\sigma^2+L_f^2)^T} +4 T L_F^2   \sigma^2(8 D^{T} +3T)$.
\section{Proof of Theorem~\ref{T1+}}
Due to the smoothness of function $F$, by setting $\eta \leq \frac{1}{2L_F}$, we have proved that $\E\left[ \Norm{\z^k}^2\right] \leq 2(L_f^2 +\sigma^2)^{T}$ and
\begin{align*}
       &\E\left[\frac{1}{K}\sum_{k=1}^K \Norm{\nabla F(\x^k)}^2 \right] \\
  \leq & \frac{2\E\left[F(\x^1)-F(\x^{K+1})\right]}{\eta K} + \frac{1}{K}\sum_{k=1}^K \E\left[\Norm{\nabla F(\x^k) - \z^k}^2\right]  - \frac{1}{2K}\E\left[\sum_{k=1}^K\Norm{\z^k}^2\right] \\
  \leq & \frac{2\E\left[F(\x^1)-F(\x^{K+1})\right]}{\eta K}   - \frac{1}{2K}\E\left[\sum_{k=1}^K\Norm{\z^k}^2\right]+D_0\alpha + 4L_F^2\left(1+ T  D^T\right) \frac{\eta^2}{\alpha^2K}\sum_{k=1}^{K+1}\E\left[\Norm{ \z^{k}}^2\right]\\
    \leq & \mathcal{O}\left( \frac{1}{\eta K} + \alpha +\frac{\eta^2}{\alpha^2 K} \right)
    \leq \mathcal{O}\left(\frac{1}{\sqrt{K}} \right),
    \end{align*}
where we set $\alpha =\mathcal{O}(K^{-1/2})$ and learning rate $\eta \leq  \frac{\alpha}{2L_F \sqrt{2(1+TD^T)}}=\mathcal{O}(K^{-1/2})$. We complete the proof of the theorem by noticing $\E\left[\frac{1}{K}\sum_{k=1}^K \Norm{\nabla F(\x^k)} \right] \leq \sqrt{\E\left[\frac{1}{K}\sum_{k=1}^K \Norm{\nabla F(\x^k)}^2 \right]}$.

\section{Dependence Issues in Previous Literature}
The main problem in the analysis of~\citet{doi:10.1137/21M1406222} lies in the variable dependence in the analysis. For example, let us examine their proof for Proposition 3.4~(refer to Proposition 3.1 if using the arXiv version). In the proof, they use
\begin{align*}
&\E \left[\left\|f_{i}\left(w_{i-1}^{k+1}\right)-w_{i}^{k+1}\right\|\left\|w_{i}^{k+1}-w_{i}^{k}\right\|^{2}| \mathcal{F}_{k}\right] \\ 
\leq & \left(\mathbb{E}\left[\left\|f_{i}\left(w_{i-1}^{k+1}\right)-w_{i}^{k+1}\right\|^{2} | \mathcal{F}_{k}\right]\right)^{\frac{1}{2}} \mathbb{E}\left[\left\|w_{i}^{k+1}-w_{i}^{k}\right\|^{2} | \mathcal{F}_{k}\right].
\end{align*}
This inequality is not true since we cannot ensure $\E \left[AB\right] \leq \E[A] \cdot  \E[B]$ when $A$ and $B$ are not independent~(e.g., when $A=B=X$ and $\E[X]=0$, it is obvious that we can not ensure $\E [X^2]\leq \E [X] \E [X]=0$ ). Note that $w_{i}^{k+1}$ appears in both terms. Thus, we cannot ensure 
\begin{align*}
&\E \left[\left\|f_{i}\left(w_{i-1}^{k+1}\right)-w_{i}^{k+1}\right\|\left\|w_{i}^{k+1}-w_{i}^{k}\right\|^{2}| \mathcal{F}_{k}\right]  \\
\leq &\E \left[\left\|f_{i}\left(w_{i-1}^{k+1}\right)-w_{i}^{k+1}| \mathcal{F}_{k}\right\|\right]\E \left[\left\|w_{i}^{k+1}-w_{i}^{k}\right\|^{2}| \mathcal{F}_{k}\right]  \\
\leq & \left(\mathbb{E}\left[\left\|f_{i}\left(w_{i-1}^{k+1}\right)-w_{i}^{k+1}\right\|^{2} | \mathcal{F}_{k}\right]\right)^{\frac{1}{2}} \mathbb{E}\left[\left\|w_{i}^{k+1}-w_{i}^{k}\right\|^{2} | \mathcal{F}_{k}\right].
\end{align*}
As a result, the analysis is invalid without this inequality.

\section{Additional Experimental Results}\label{addExp}

\begin{figure*}[!ht]
\begin{center}
\includegraphics[width=\rata\textwidth]{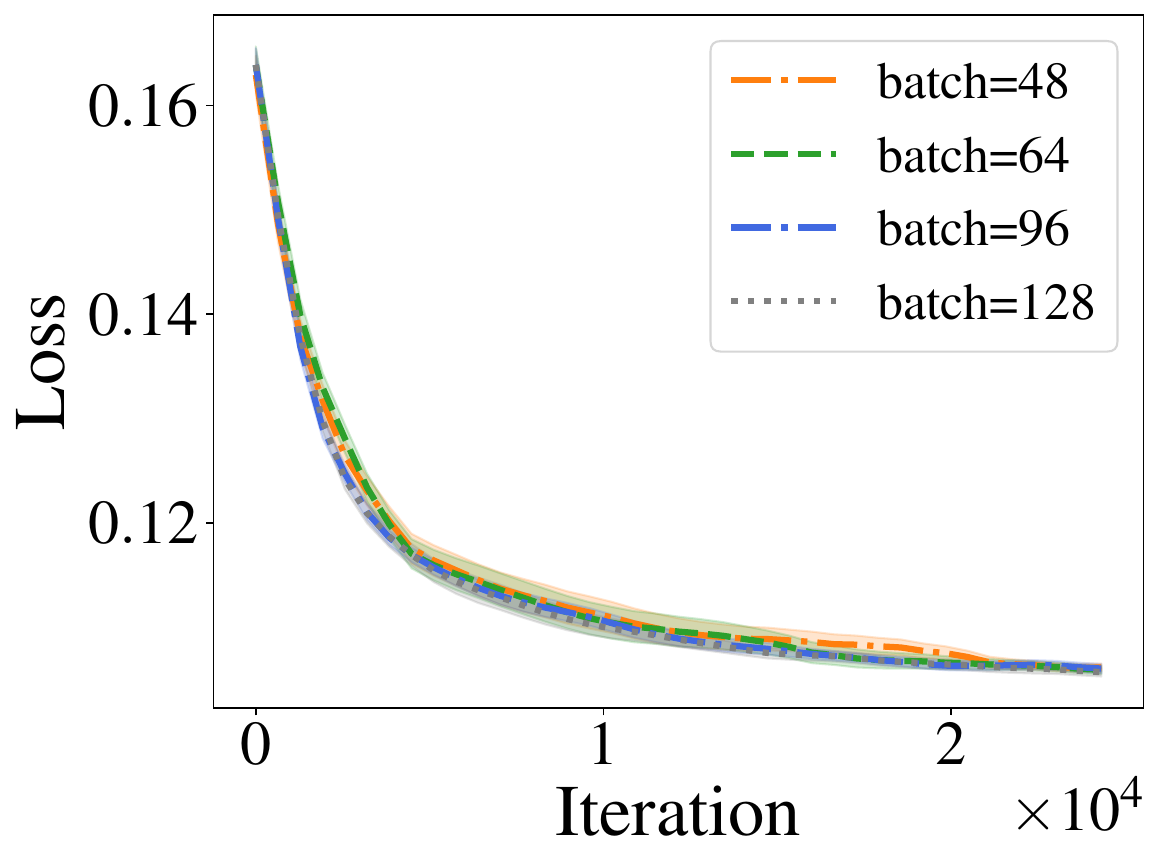}
\includegraphics[width=\rata\textwidth]{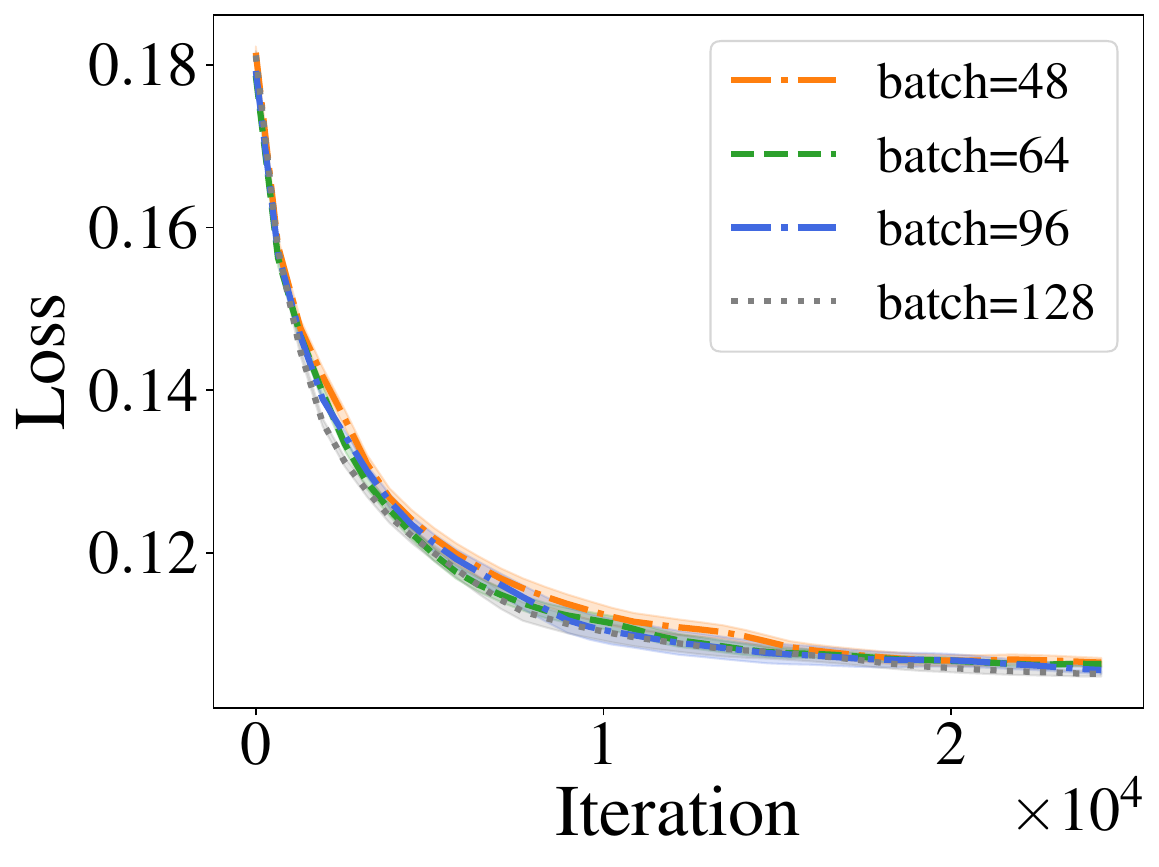}
\includegraphics[width=\rata \textwidth]{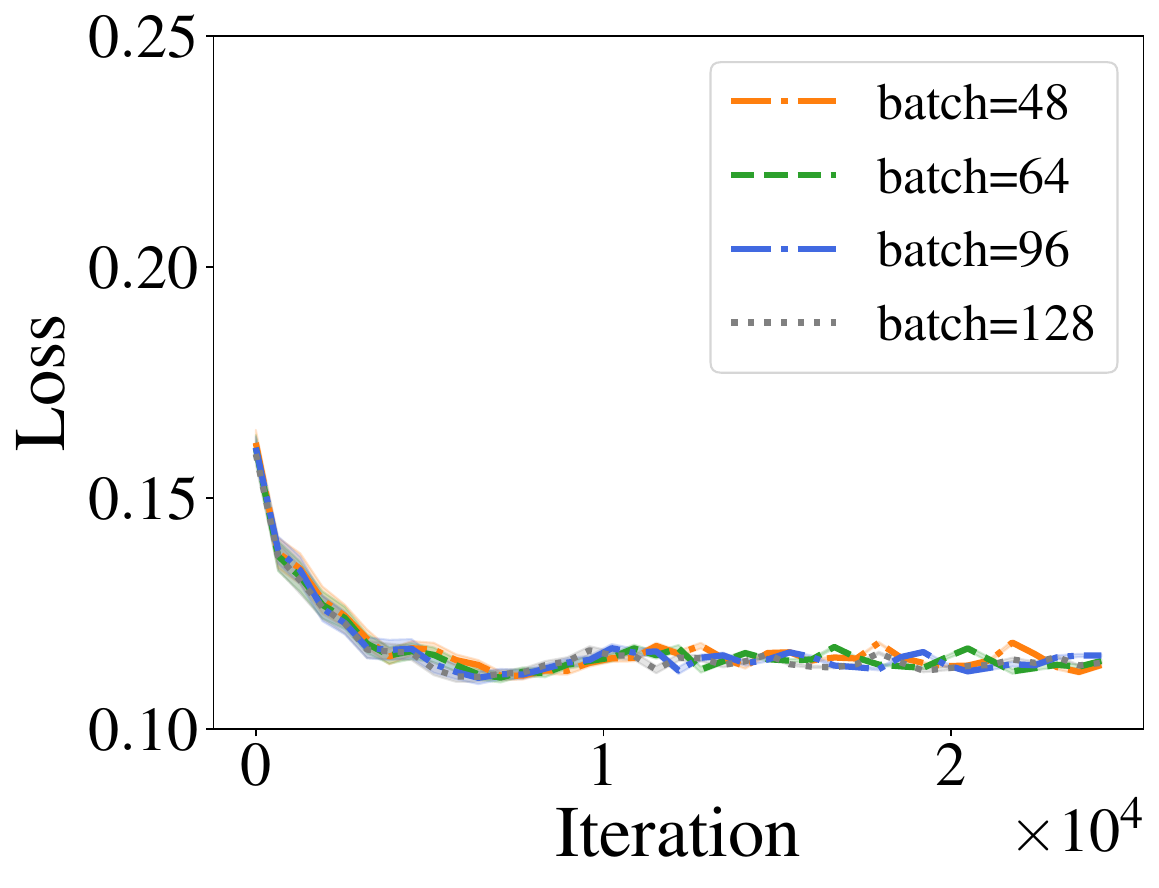}\\
\subfigure[Industry-10]{\includegraphics[width=\rata\textwidth]{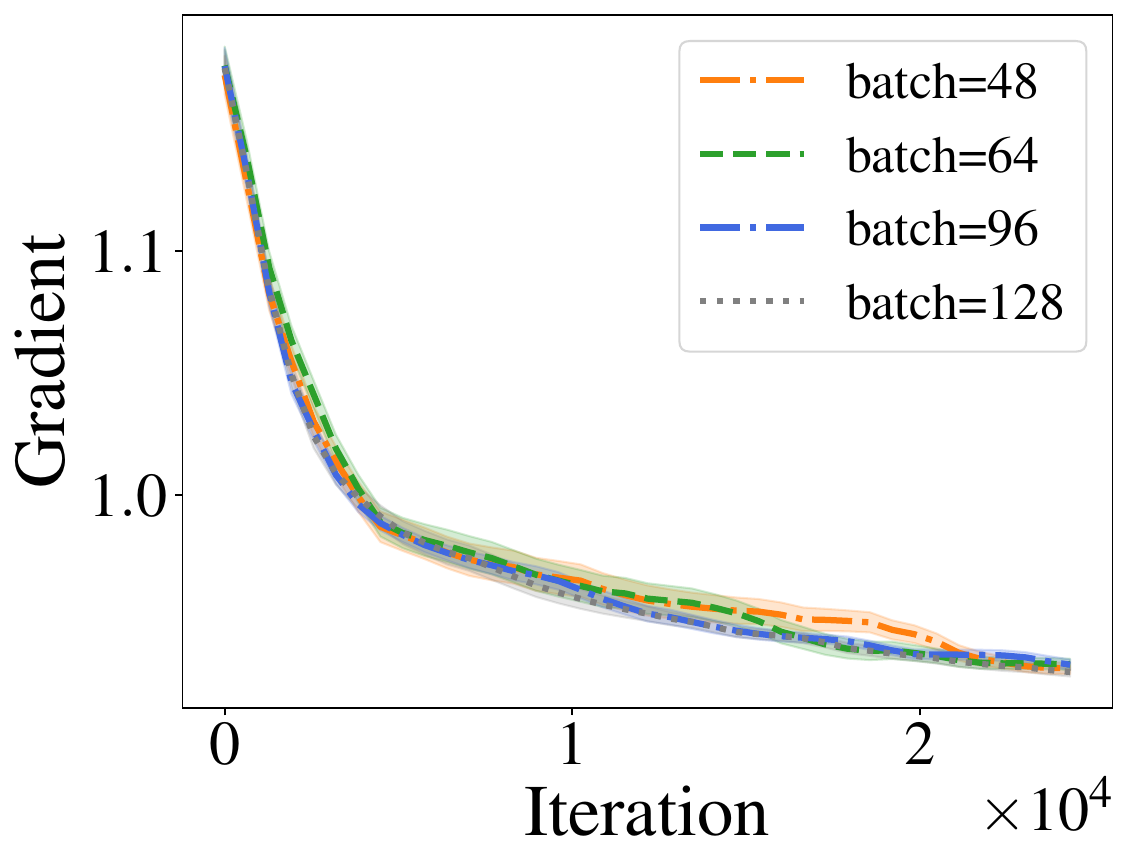}}
\subfigure[Industry-12]{\includegraphics[width=\rata\textwidth]{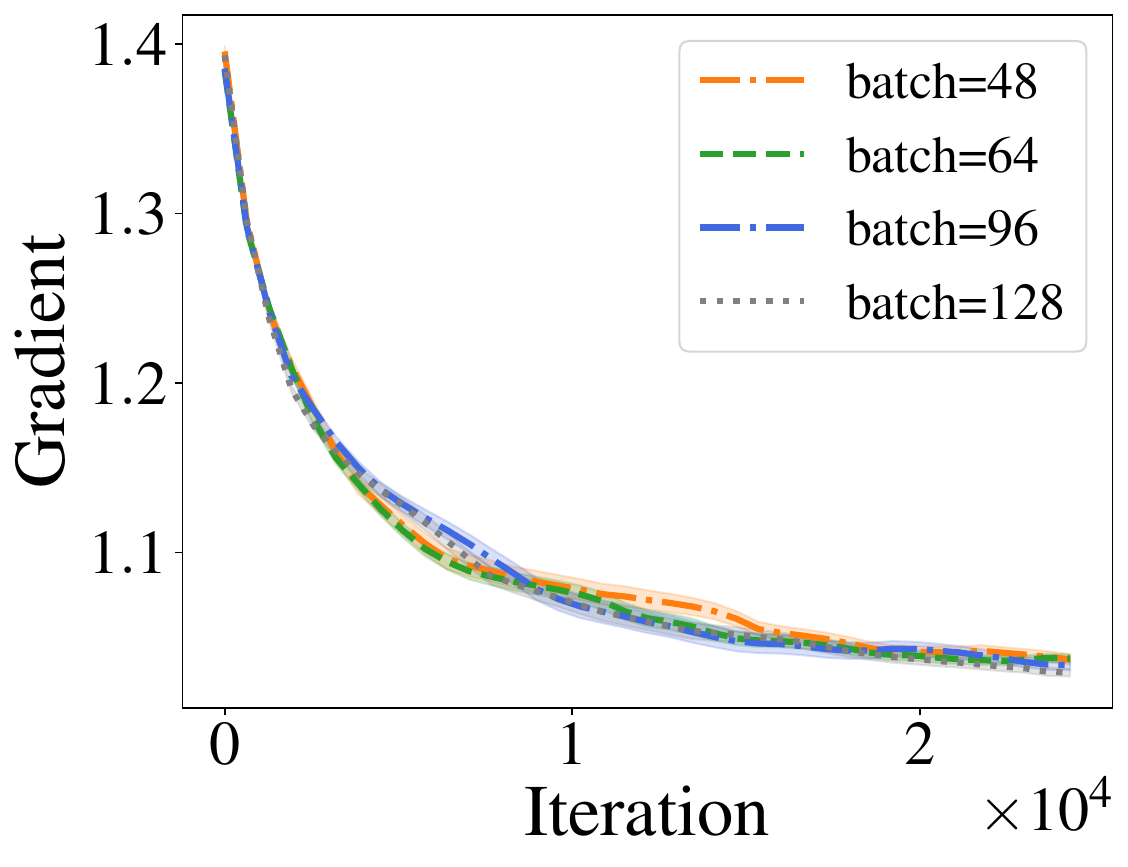}}
\subfigure[Industry-17]{\includegraphics[width=\rata\textwidth]{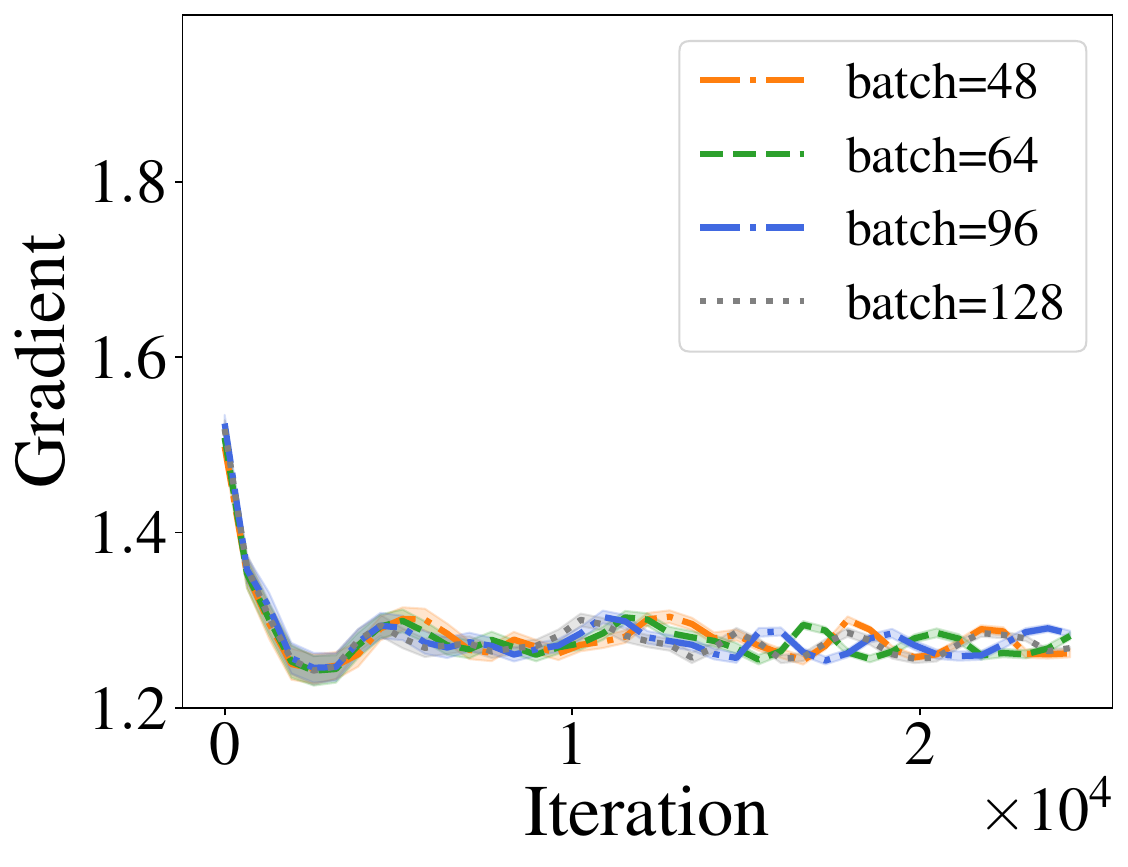}}
\caption{Results for Risk-Averse Portfolio Optimization with different batch sizes.}
		\label{fig:01}
	\end{center}
\end{figure*}

\begin{figure*}[!ht]
  \begin{center}
\includegraphics[width=\rata\textwidth]{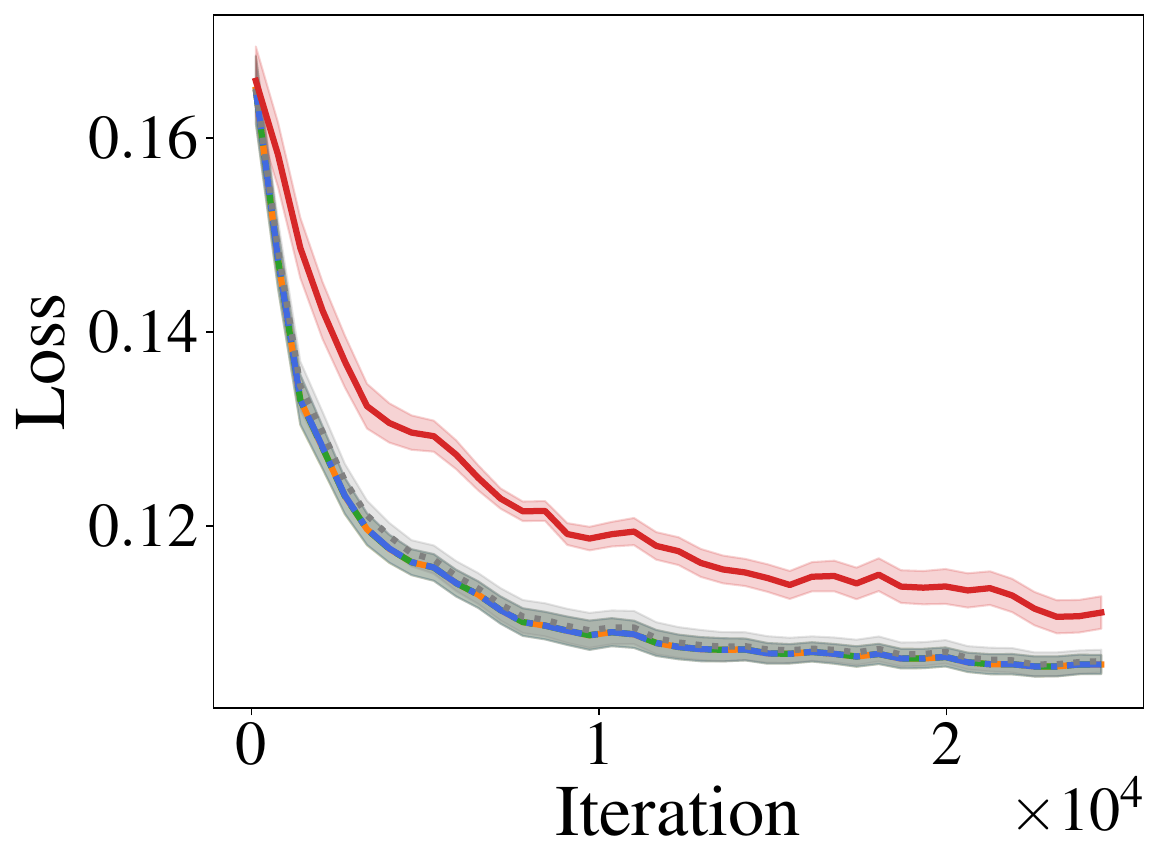}
\includegraphics[width=\rata\textwidth]{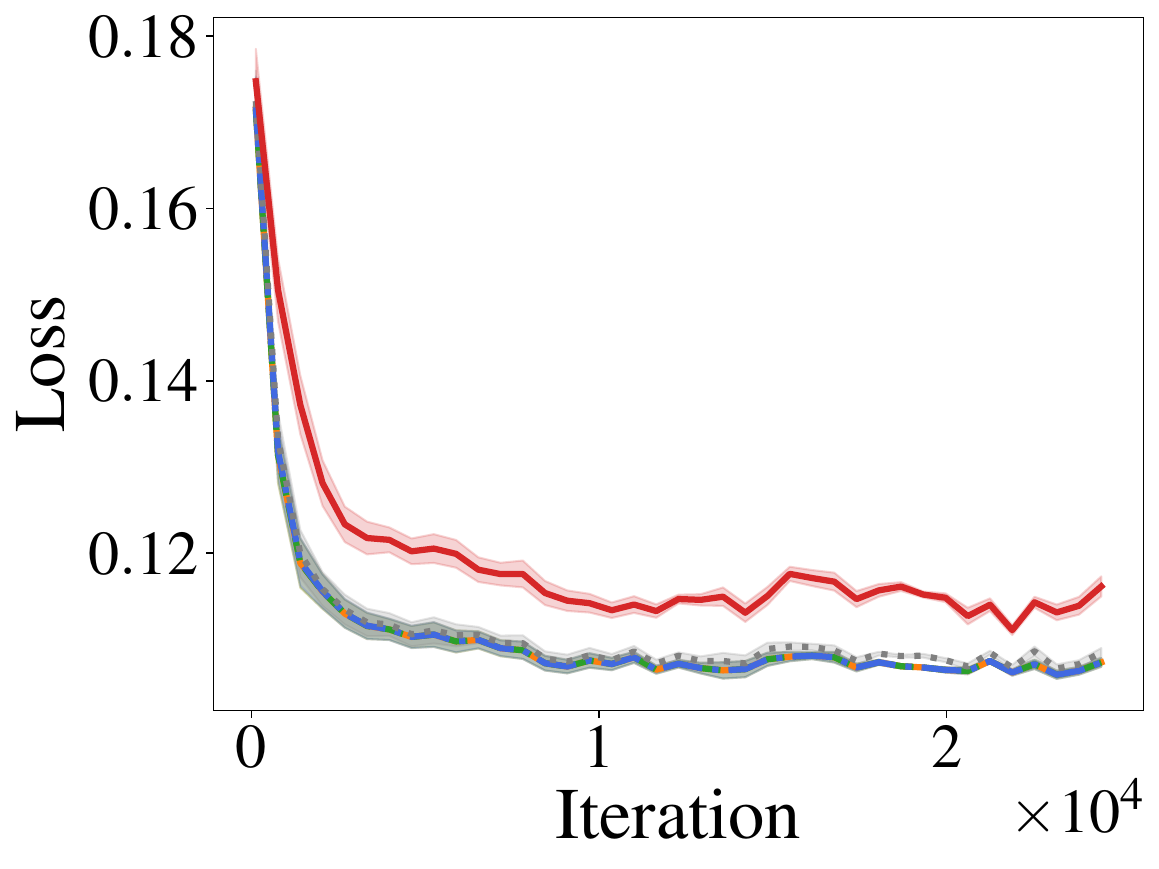}
\includegraphics[width=\rata\textwidth]{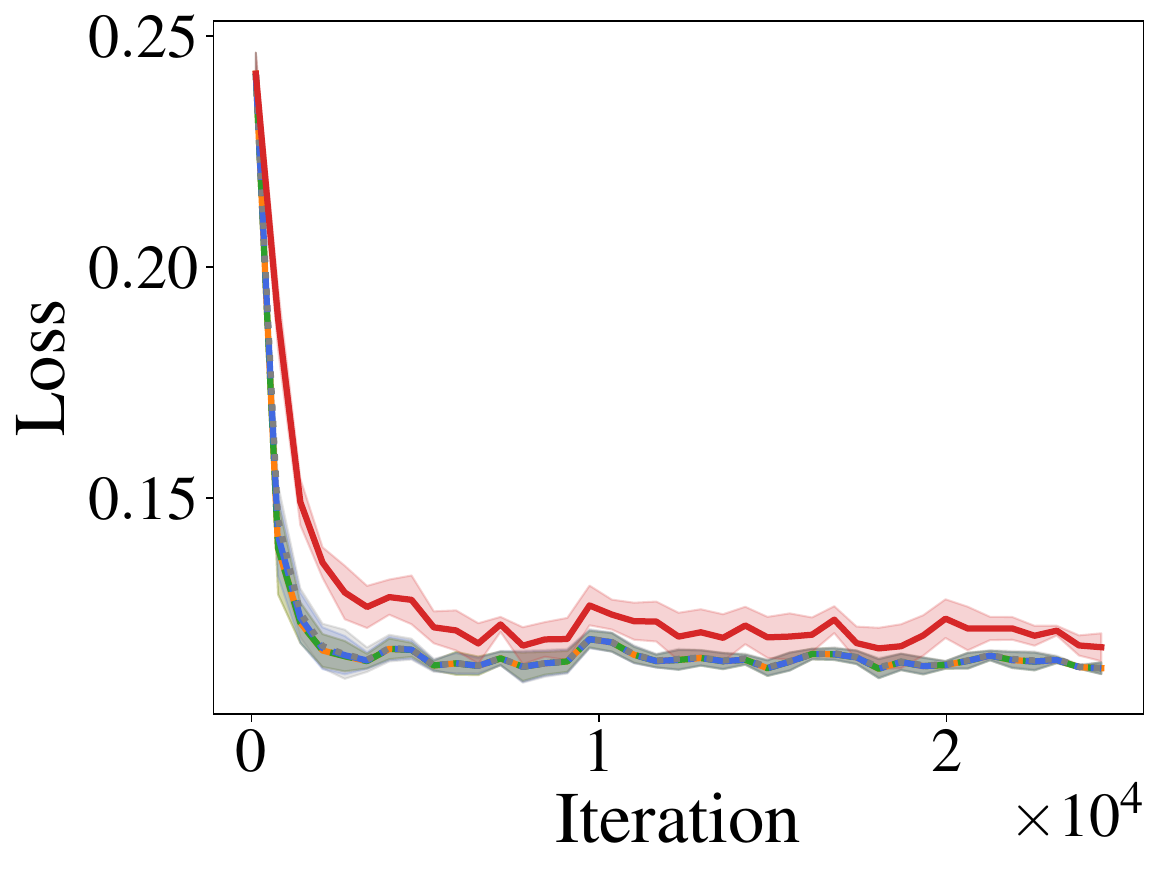}\\
\subfigure[Industry-10]{\includegraphics[width=\rata\textwidth]{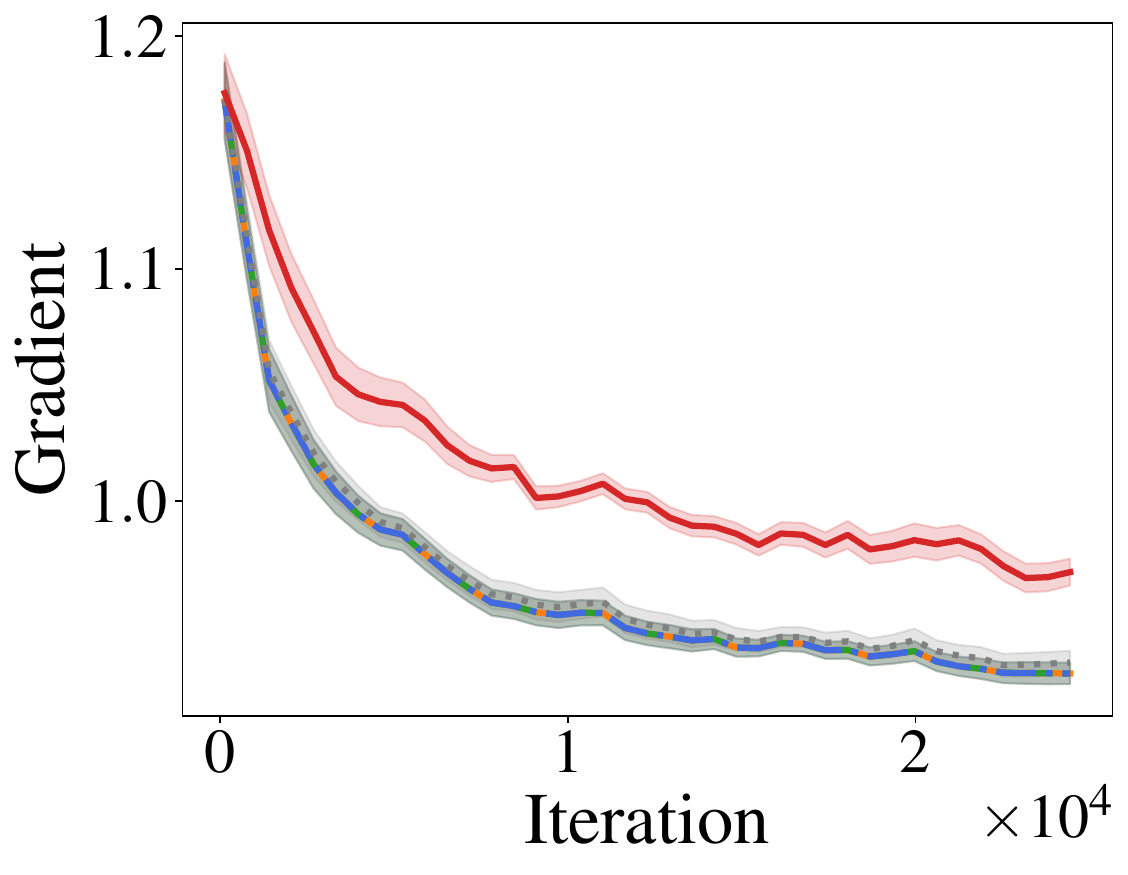}}
\subfigure[Industry-12]{\includegraphics[width=\rata\textwidth]{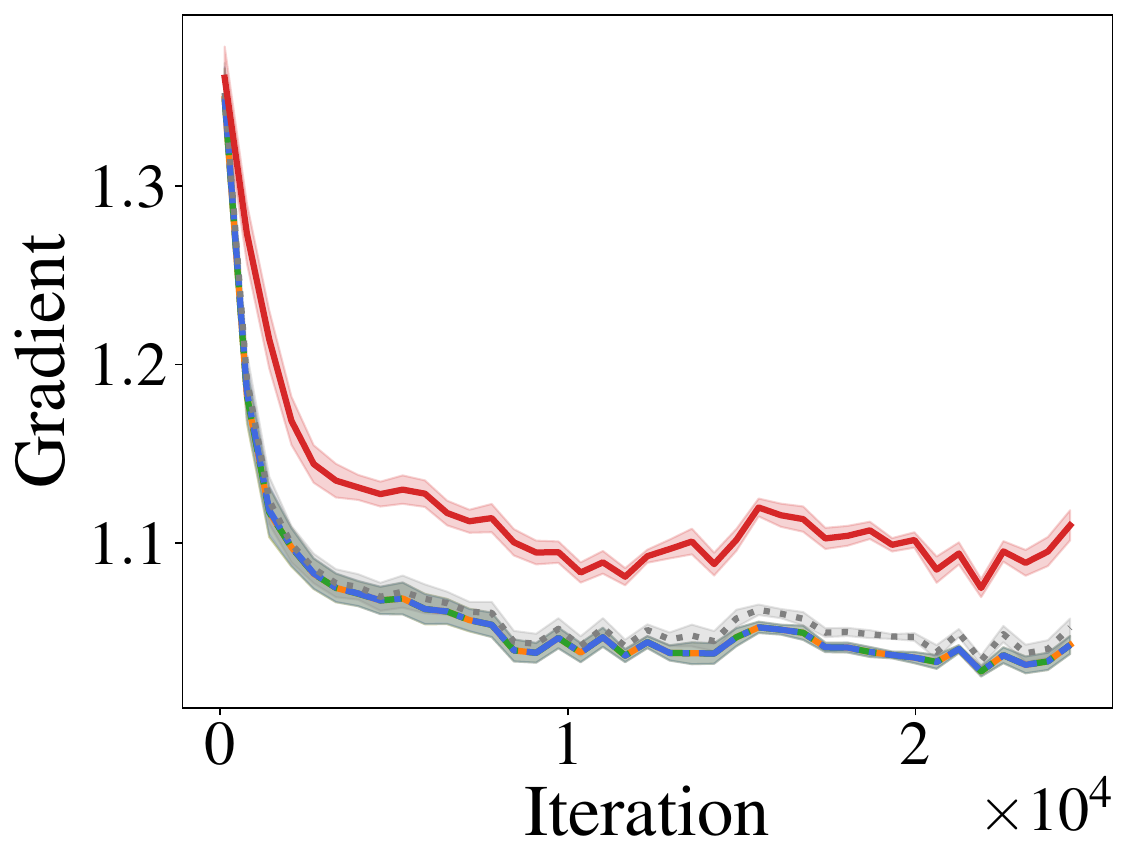}}
\subfigure[Industry-17]{\includegraphics[width=\rata\textwidth]{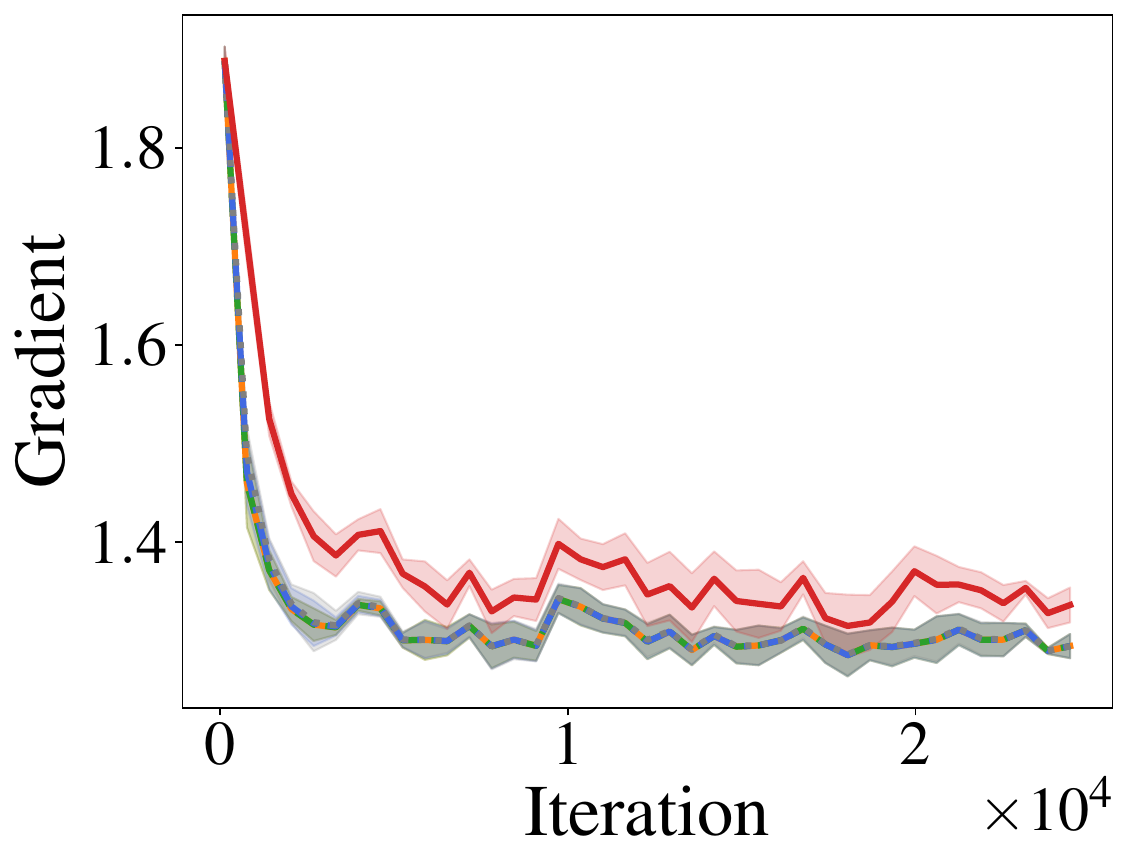}}
\includegraphics[width=0.75\textwidth]{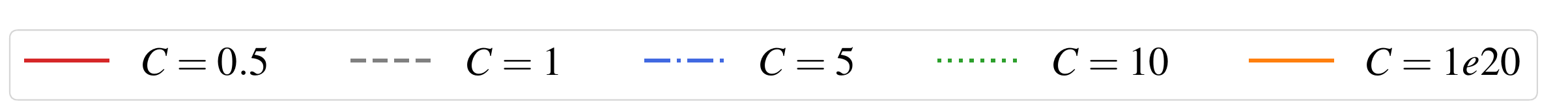} 
\caption{Sensitivity to the clipping threshold $C$ on the portfolio task.}
    \label{fig:2}
  \end{center}
\end{figure*}

\begin{figure*}[!ht]
  \begin{center}
\includegraphics[width=\rata\textwidth]{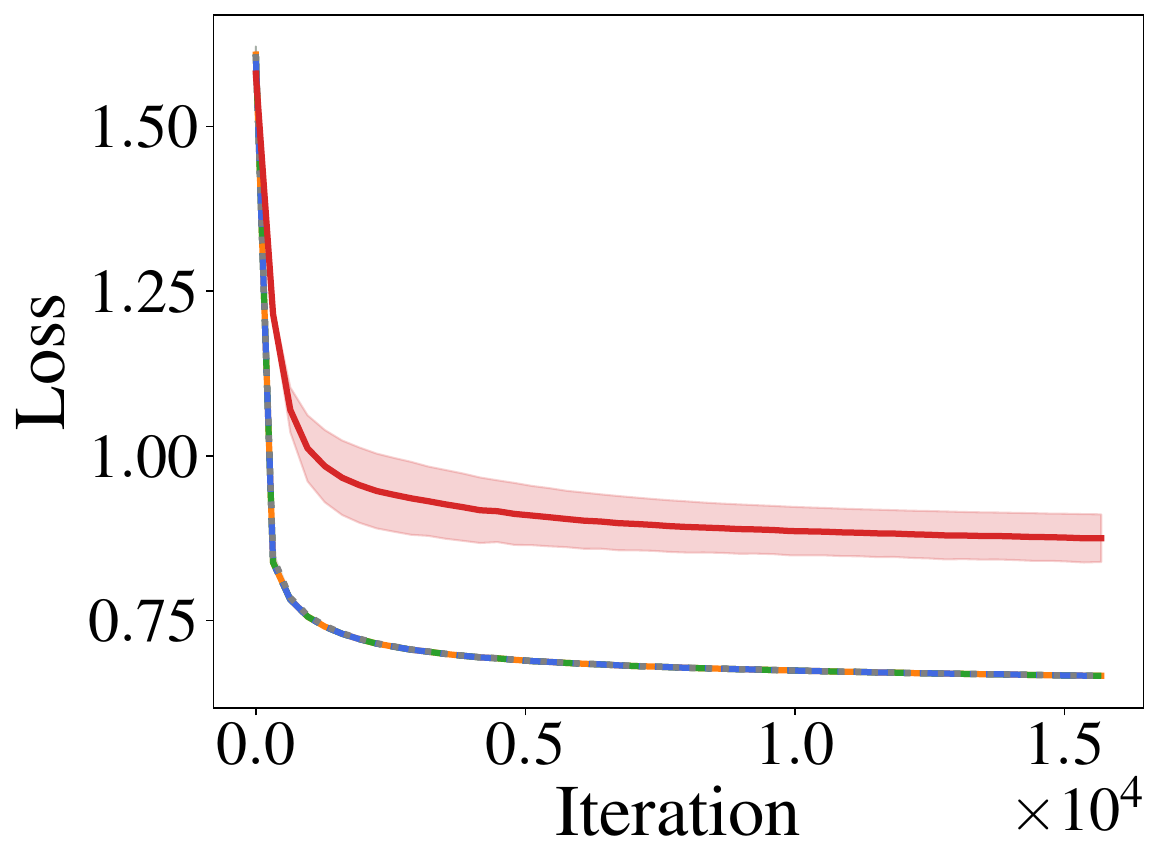}   
\includegraphics[width=\rata\textwidth]{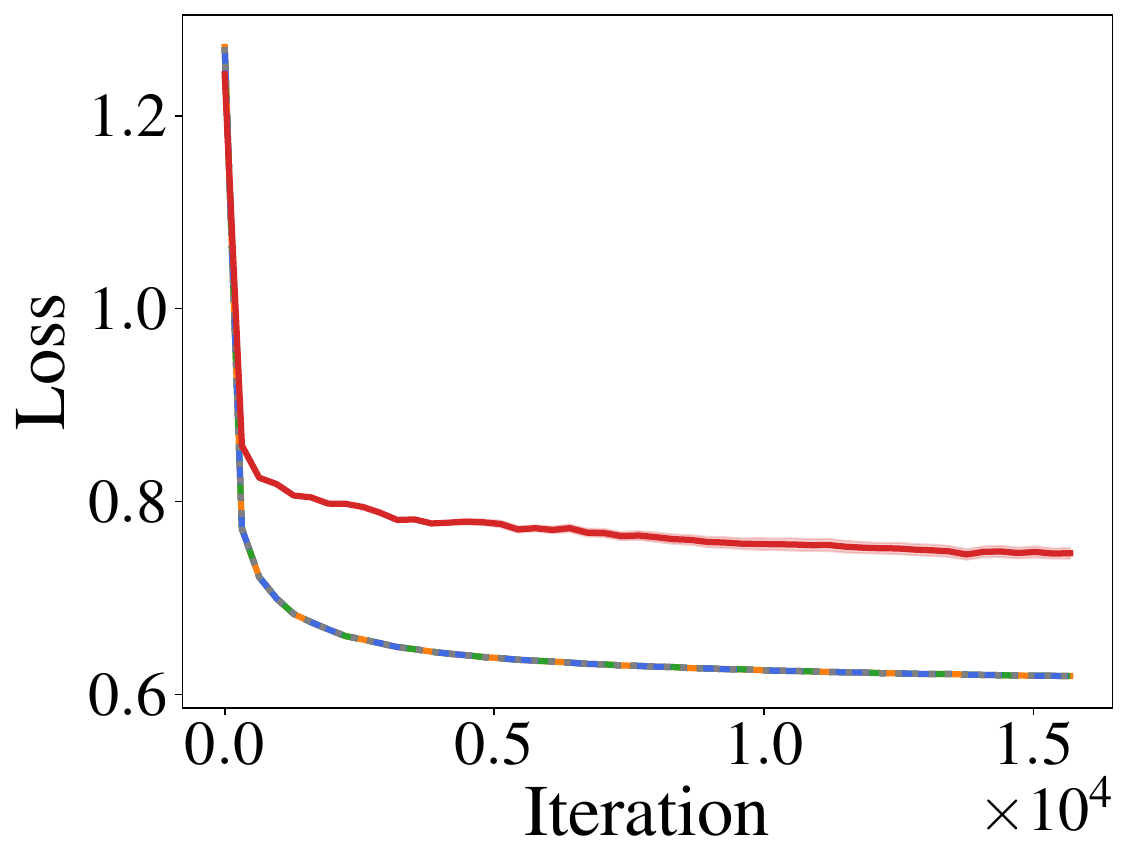}   
\includegraphics[width=\rata\textwidth]{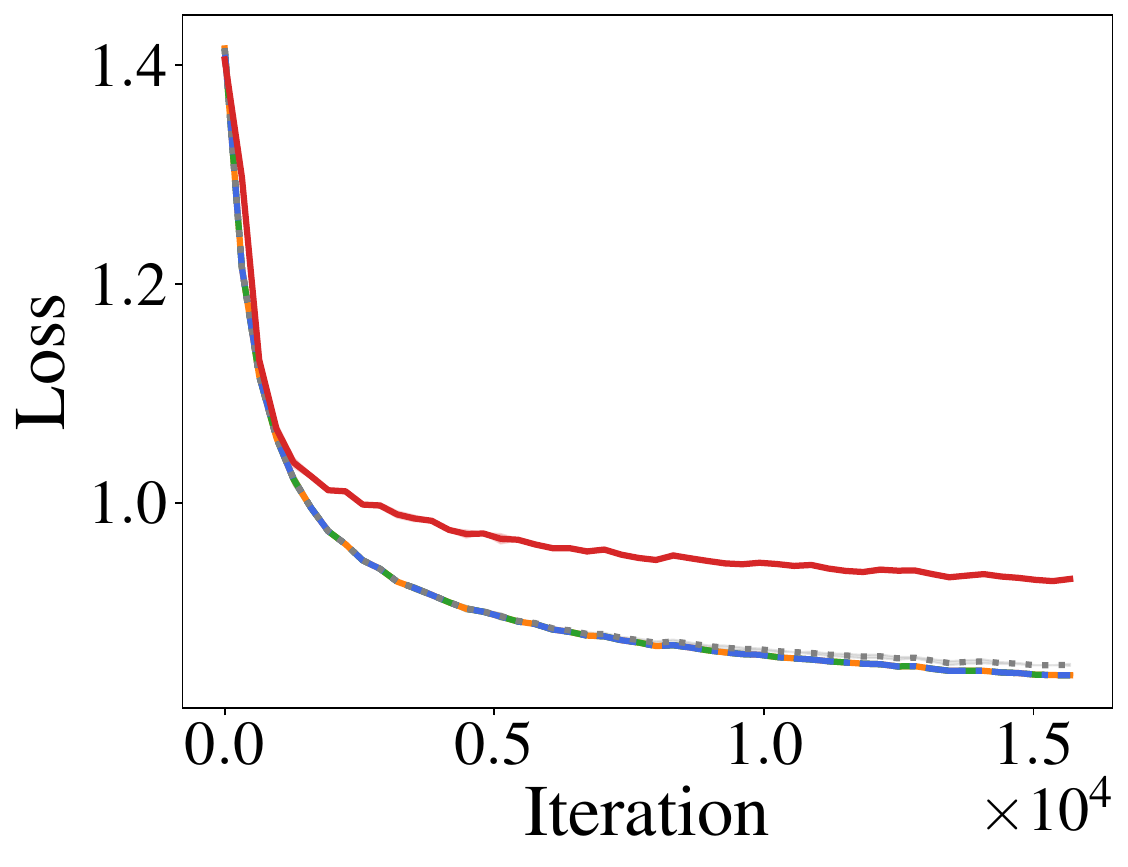}    \\
\subfigure[Australian]{\includegraphics[width=\rata\textwidth]{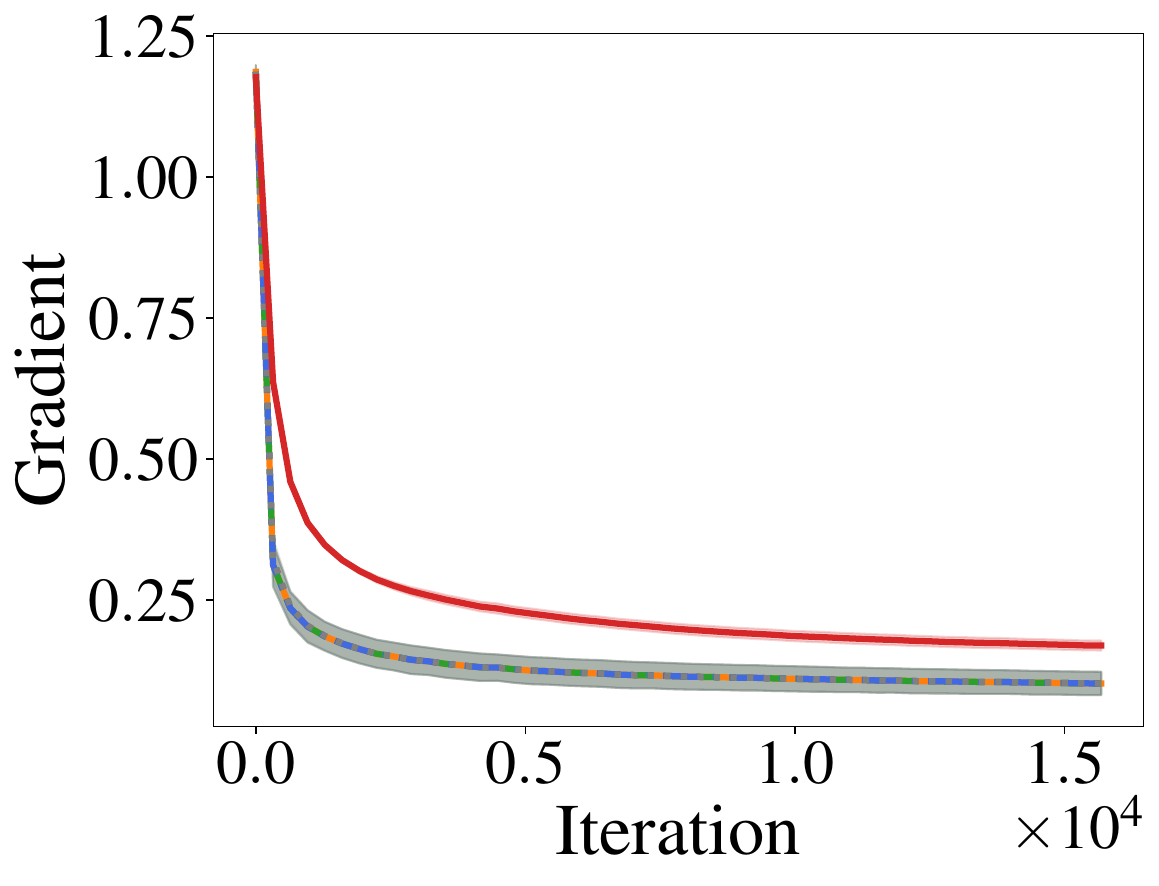}}
\subfigure[Breast-cancer]
{\includegraphics[width=\rata\textwidth]{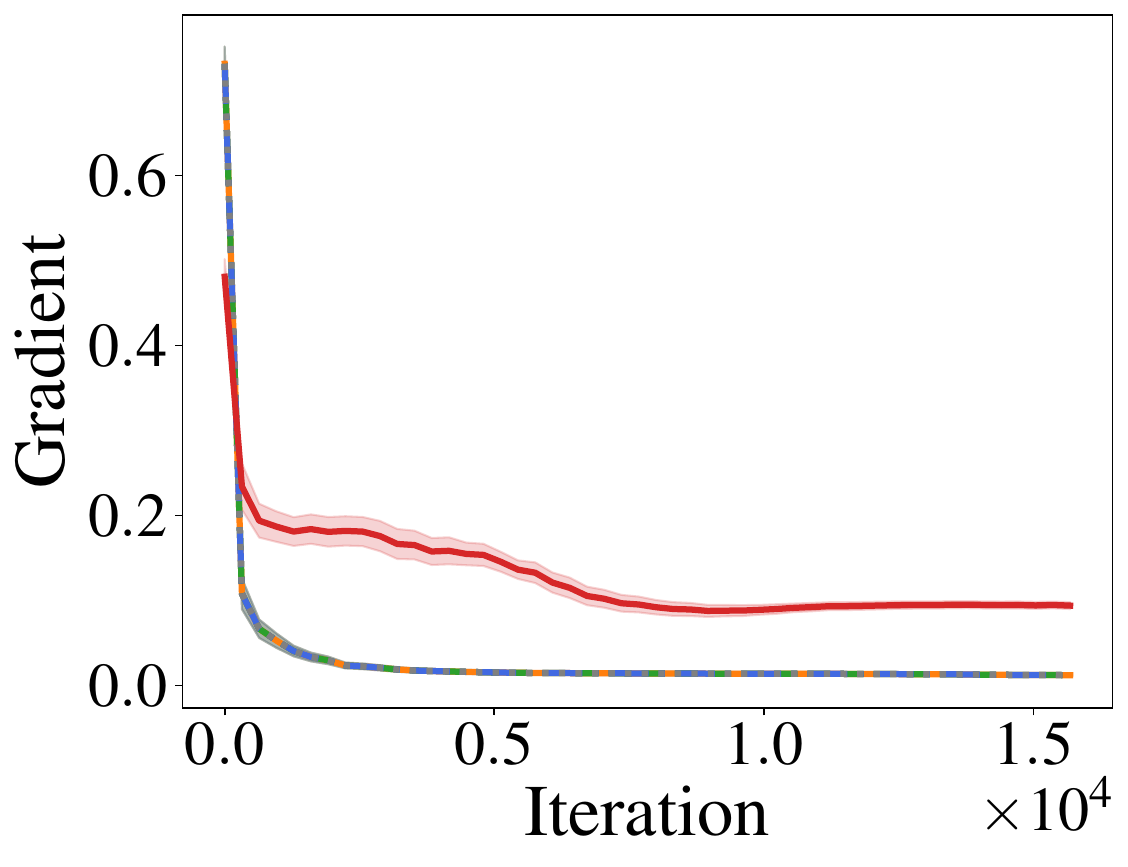}} 
\subfigure[Svmguide1]{\includegraphics[width=\rata\textwidth]{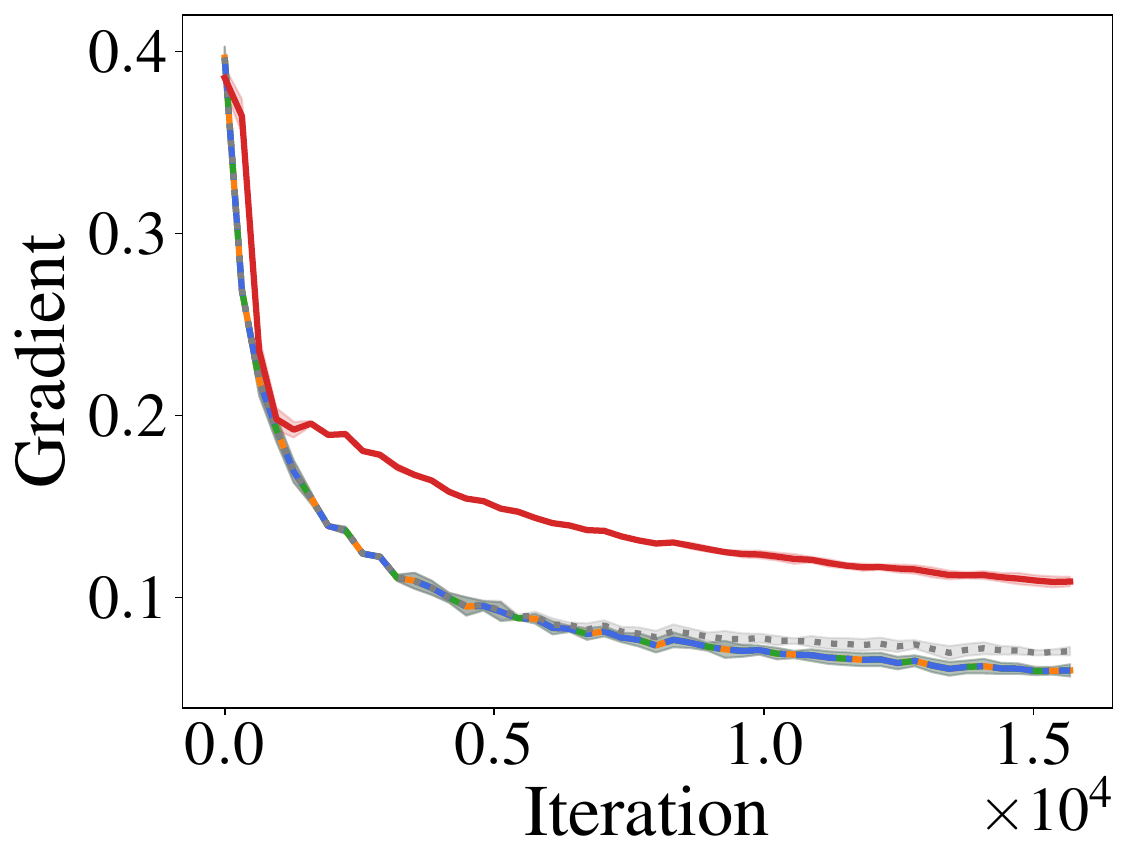}} 
\includegraphics[width=0.75\textwidth]{./legend2.pdf}
\caption{Sensitivity to the clipping threshold $C$ on hierarchical TERM.}
    \label{fig:4}
  \end{center}
\end{figure*}
\end{document}